
\input amstex
\input amsppt.sty
\magnification 1200
\NoBlackBoxes
\hcorrection{0.25in}

\topmatter
\title Log Sarkisov Program\endtitle
\author Andrea Bruno and Kenji Matsuki\endauthor
\endtopmatter

\document
\baselineskip 14pt

$$\bold{\text{Contents}}$$

\S 0. Introduction.

\S 1. Flowchart for Sarkisov Program.

\S 2. Termination of Flowchart.

\S 3. Log Sarkisov Program with KLT Singularities.

\S 4. Log Sarkisov Program with WKLT Singularities.

\vskip.2in

$$\footnote""{Research at MSRI is supported in part by NSF grant no.
DMS-9022140}$$

The purpose of this paper is two-fold.  The first is to give a
tutorial introduction to the so called Sarkisov program, a 3-dimensional
generalization of Castelnuovo-N\"other Theorem ``untwisting" birational maps
between Mori fiber spaces, which was recently established by Corti[4].  We
should emphasize that though the general features were understood
(cf.Matsuki[17]) after
Reid[29] explained the original ideas of V.G. Sarkisov in a
substantially laundered form, it is only after Corti[4]
that we are beginning to understand the details of the mechanism.  Here we will
present a flowchart to visualize how the Sarkisov program works and also
slightly simplify the proof of termination
after the ingenious argument of Corti[4]: we prove there is no infinite loop in
the program just observing that the Sarkisov degree decreases strictly after
each
untwisting and it cannot decrease infinitely many times using the boundedness
of
${\Bbb Q}$-Fano $d$-folds $d \leq 3$ together with the ascending chain
condition
$S_3(\text{Local})$ of Alexeev[1] (cf.Shokurov[31]Koll\'ar et al[16]).  Our
argument also makes it explicit that the Sarkisov program holds in arbitrary
dimension $n$ once we have Log MMP in dimension $n$, boundedness of ${\Bbb
Q}$-Fano $d$-folds for $d \leq n$ and $S_n(\text{Local})$.  The second is an
attempt to give a logarithmic generalization following the philosophy of
Iitaka,
based upon the Log MMP (established in dimension 3 by Shokurov[31]Kawamata[9]
(cf.Koll\'ar et al[16])).  The key is to understand the meaning of the
genuine Sarkisov program and set up the natural and right generalization.  The
genuine Sarkisov program untwists any birational map between two Mori fiber
spaces which are birationally equivalent.  A naive speculation that the log
Sarkisov program should untwist any birational map between two log Mori fiber
spaces which are only birationally equivalent turns out to be not natural and
simply does not work!  In order to reach the right understanding of what the
genuine Sarkisov program does we have to introduce the notion of the Sarkisov
relation:Mori fiber spaces are Sarkisov related iff they are the end results of
the $K$-MMP starting from an appropriate nonsingular projective variety.  The
genuine Sarkisov program untwists a birational map between two Mori fiber
spaces
which are Sarkisov related, by factorizing it into ``links" among intermediate
Mori fiber spaces, all of which including the original two are Sarkisov
related.  Log Mori fiber spaces are said to be Sarkisov related iff they are
the
end results of the $K + B$-MMP starting from one log pair consisting of a
nonsingular variety and an S.N.C. divisor as a boundary.  The log Sarkisov
program should be the one to untwist a birational map between two log Mori
fiber
spaces which are Sarkisov related, by factorizing it into ``links" among
intermediate log Mori fiber spaces, all of which including the original two are
Sarkisov related.  Once this is understood, the log Sarkisov program works
almost parallel to the genuine Sarkisov program in the case of kawamata log
terminal singularities in arbitrary dimension, except for the verification of
termination.  The boundedness of a certain class of log
${\Bbb Q}$-Fano varieties becomes crucial for our argument in showing
termination, just as the result of Kawamata[8] on the boundedness of ${\Bbb
Q}$-Fano 3-folds was crucial in showing termination of the genuine Sarkisov
program in dimension 3.  While in dimension 2 we establish termination thanks
to
a result of Nukulin[26]Alexeev[2], in dimension 3 we have to use a conjecture
by
Borisov[3] at the last stage of the proof of termination.  We also establish
the
log Sarkisov program in the case of weakly kawamata log terminal singularities
in
dimension 2 including its termination.  But in dimension 3, the Sarkisov
relation becomes quite subtle for wklt singularities and we can only discuss
problems toward establishing the program including not only its temination but
the general mechanism itself.

We also remark that the birational transformations among the various moduli
spaces studied by
M. Thaddeus and others can be put into the general frame work of the (log)
Sarkisov program.

We would like to thank A. Corti, who allowed us to present many of his
ideas here in logarithmic form.  Most of the arguments here are taken from his
paper Corti[4]
and repeated for the sake of understanding of the reader.  The original ideas
are due entirely to V.G. Sarkisov and our indebtedness to M. Reid toward
understanding them is as clear as his paper Reid[29].  The conversations with
J. Koll\'ar and J. McKernan were very helpful and critical.  We would like to
thank S. Mori, who gave us warm encouragement throughout.

\vskip.2in

\subhead{\S 0. Introduction}\endsubhead

\vskip.1in

{}From the view point of Minimal Model Program (the so-called Mori's program),
our basic strategy to understand the birational geometry of higher
dimensional algebraic varieties (established in dimension 3 and conjectural
for higher dimension) is divided into the following 3 steps:

\ \ 1.  Find ``good representatives" among varieties with a given function
field through MMP:  We take a nonsingular projective variety $X$ by Hironaka's
resolution of singularities with a given function field.  Input $X$ into the
black box called Minimal Model Program (abbreviated MMP), which produces a good
representative, i.e., a minimal model or a Mori fiber space as an output,
depending whether $X$ is non-uniruled or uniruled.

$$
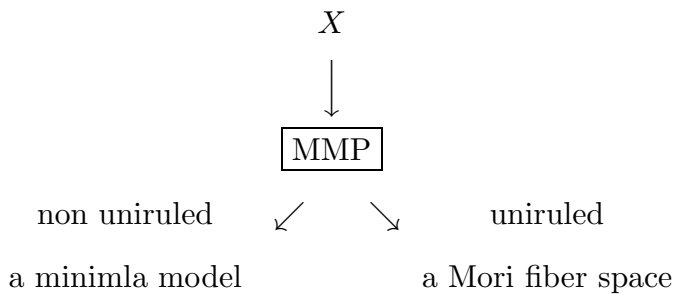
\CD
@.X@.\\
@.@VVV @.\\
@.\boxed{\text{MMP}}@.\\
\text{non\ uniruled}@.\ \ \ \swarrow\ \ \ \ \ \ \searrow\ \ @.\text{uniruled}\\
\text{a\ minimla\ model}@.@.\text{a\ Mori\ fiber\ space}\\
\endCD
$$

\ \ 2. Study the properties of ``good representatives": The most important
property is the Dichotomy, which says that the uniruledness should
characterize the variety with Kodaira dimension $- \infty$, i.e. $\kappa = -
\infty$ or $\kappa \geq 0$ depending upon whether $X$ is uniruled or not.  In
dimension 3, this is a theorem.  For 3-folds $X$ with $\kappa(X)\geq 0$, the
Abundance Theorem of Kawamata-Miyaoka further claims that a minimal model
$X_{min}$ has a base point free pluri-canonical system which induces the
canonical morphism $\Phi_{|mK_{X_{min}}|}:X_{min} \rightarrow X_{can}$,
crucial for the understanding of the global structure of $X$ and its moduli.
For 3-folds with $\kappa = - \infty$, a theorem of Miyaoka-Mori says that a
Mori
fiber space $X_{mori}$ is covered by rational curves intersecting
$K_{X_{mori}}$ negatively.

\ \ 3. Study the relation among ``good representatives": For 3-folds with
$\kappa \geq 0$, the basic relation among good representatives is that minimal
models in a given birational equivalence class are connected by a sequence of
flops (cf. Reid[28]Kawamata[5]\linebreak
Koll\'ar[15].  See also Matsuki[18] for a finer
description of their relation).  It is the Sarkisov program, the main theme of
our paper, which describes the relation among good representatives for 3-folds
with $\kappa =
- \infty$, i.e., Mori
fiber spaces in a given birational equivalence class.

\vskip.1in

In dimension 2, i.e., in the case of classical birational geometry of
surfaces, the meaning of these 3 steps is rather straightforward.

\ \ 1. Starting from a nonsingular projective surface, we keep contracting
(-1)-curves (MMP in dimension 2) until we get either a surface $X_{min}$ with
the canonical divisor $K_{X_{min}}$ being nef or a ruled surface $X_{mori}$
over a curve (or ${\Bbb P}^2$ over a point).

\ \ 2. When $\kappa(X) \geq 0$, the canonical morphism
$\Phi_{|mK_{X_{min}}|}:X_{min} \rightarrow X_{can}$ from a minimal model is:

$\kappa = 2$ - a biational map to a canonically polarized surface with only
rational double points

$\kappa = 1$ - an elliptic fibtration whose degeneration fibers are
studied by Kodaira[14]

$\kappa = 0$ - a trivial map to a point, where we know $X_{min}$
must be either Abelian, bielliptic, K3 or Enriques.

When $\kappa = - \infty$, the structure of a Mori fiber space is rigid and
well-understood: either a ${\Bbb P}^1$-bundle over a nonsingular curve or
${\Bbb P}^2$.

\ \ 3. A minimal model is unique in a fixed birational equivalence class for
surfaces
with $\kappa \geq 0$, while any birational map among ruled surfaces in a given
birational
equivalence class is decomposed into a sequence of elementary transformations
by
Castelnuovo-N\"other theorem.

\vskip.1in

$$\bold{Logarithmic\ Generalization}$$

The logarithmic generalization of the basic strategy following the philosophy
of Iitaka to understand the birational geometry of varieties WITH BOUNDARIES
goes along the same line:

\ \ 1. We take a pair $(X,B_X)$ where $X$ is a nonsingular projective
variety and $B_X = \Sigma b_iB_i$ is a simple normal crossing divisor with
$0 \leq b_i \leq 1$.  Input $(X,B_X)$ into the black box called Log MMP,
which produces a log minimal model $(X_{min},B_{X_{min}})$ with
$K_{X_{min}} + B_{X_{min}}$ being nef or a log Mori fiber space
$\phi:(X_{mori},B_{X_{mori}}) \rightarrow S$ with $K_{X_{mori}} +
B_{X_{mori}}$ being $\phi$-negative, depending upon whether $(X,B_X)$ is
non log uniruled or log uniruled.  (We say $(X,B_X)$ is log uniruled iff it is
covered by rational curves intersecting $K_X + B_X$ negatively.)

$$
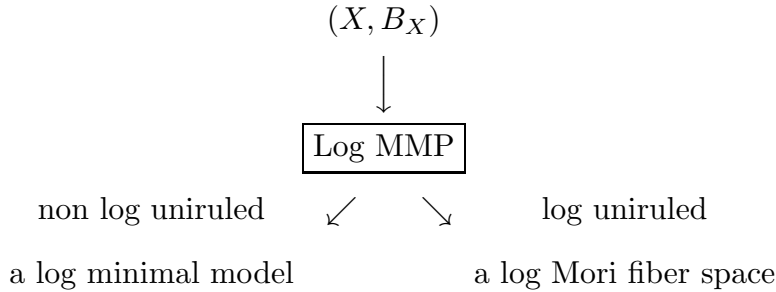
\CD
@.(X,B_X)@.\\
@.@VVV @.\\
@.\boxed{\text{Log\ MMP}}@.\\
\text{non\ log\ uniruled}@.\ \ \ \swarrow\ \ \ \ \ \ \searrow\ \ @.\text{log\
uniruled}\\
\text{a\ log\ minimal\ model}@.@.\text{a\ log\ Mori\ fiber\ space}\\
\endCD
$$

\ \ 2. Again the most important
property is the Dichotomy, which says that the log uniruledness should
characterize the varieties with log Kodaira dimension $- \infty$, i.e.
$\kappa(K_X + B_X) = - \infty$ or $\kappa(K_X + B_X) \geq 0$ depending
upon whether
$(X,B_X)$ is log uniruled or not.  In
dimension 3, this is a theorem.  For log 3-folds $(X,B_X)$ with
$\kappa(K_X + B_X) \geq 0$, the Log Abundance Theorem of
KeelMatsukiMcKernan[12] further
claims that a log minimal model $(X_{min},B_{X_{min}})$ has a base point free
pluri-log
canonical system.  For log 3-folds with $\kappa = - \infty$, a theorem of
Miyaoka-Mori applies again to imply that a log Mori fiber space
$(X_{mori},B_{X_{mori}})$ is covered by rational curves intersecting
$K_{X_{mori}} + B_{X_{mori}}$ negatively.

\ \ 3. In the genuine birational geometry, we are interested in the relation
among good representatives in a given birational equivalence class, where two
good representatives are outcomes of one appropriate nonsingular projective
variety through MMP if and only if they are birationally equivalent.  This is
not the case with logarithmic birational geometry.  We say (two
or more) good representatives are Sarkisov related iff they are outcomes
through Log MMP of one appropriate log pair consisting of a nonsingular
projective variety and a S.N.C. divisor as a boundary.  Then log minimal models
which are Sarkisov related are connected by a sequence of log flops (cf.
Koll\'ar[15]Koll\'ar et al[16].  See \S 4 for a detalied discussion.)  The log
Sarkisov program should be the one to untwist a birational map between two log
Mori fiber spaces which are Sarkisov related.

\newpage

\subhead{\S 1. Flowchart for Sarkisov Program}\endsubhead

\vskip.1in

In this section, we review the (genuine) Sarkisov program after
Corti[4]\linebreak
(cf.Sarkisov[30]Reid[29]Matsuki[17]) with some simplifications and
present a flowchart to visualize how it works.  The aim of this section is
mostly tutorial aside from simultaneously preparing the notations for the
logarithmic case which goes almost parallel to the genuine case after the
introduction of
the Sarkisov relation (See \S 3.).  We refer the reader to
KawamataMatsudaMatsuki[11] for the general features of MMP and to Koll\'ar et
al[16] for those of Log MMP.

The Sarkisov program, in short, is an
algorithm to factorize a birational map between two Mori fiber spaces, i.e.,
two different end results of MMP of one appropriate nonsingular projective
variety, when the Kodaira dimension is $- \infty$.

\proclaim{Definition 1.1} A Mori fiber space $\phi:X \rightarrow S$ is the
contraction of an extremal ray with respect to $K_X$ from a normal projective
variety with only ${\Bbb Q}$-factorial terminal singularities onto a variety
$S$,
i.e., $\phi$ is a morphism from a normal projective
variety with only ${\Bbb Q}$-factorial terminal singularitieswith connected
fibers onto a normal variety $S$ with
$\dim S <
\dim X$ s.t. $\rho(X/S) = 1$ and $-K_X$ is $\phi$-ample.
\endproclaim

\proclaim{Theorem 1.2 (Sarkisov Program in dimension 3)
(cf.Sarkisov[30]Reid[29]\linebreak
Corti[4])} A birational map
$$
\CD
X @.\overset\Phi\to\dashrightarrow @.X' \\
@V\phi VV @.@VV\phi' V \\
S @.@. S' \\
\endCD
$$
between two Mori fiber spaces in dimension 3 is a composite of the following 4
types of links.

Links of type (I)

$$
\matrix
& & Z & \dashrightarrow & X_1 \\
& \swarrow & & & \\
X & & & & \downarrow \\
\downarrow & & & & \\
S & & \leftarrow & & S_1 \\
\endmatrix
$$
where $Z \rightarrow X$ is a $K$-negative extremal divisorial contraction, $Z'
\dashrightarrow X_1$ a sequence of log flips with respect an appropriate log
pair and
$\rho(S_1/S) = 1$.

\vskip.1in

Links of type (II)

$$
\matrix
& & Z & \dashrightarrow & Z' & & \\
& \swarrow & & & & \searrow & \\
X & & & & & & X_1 \\
\downarrow & & & & & & \downarrow \\
S & & & \overset\sim\to\leftarrow & & & S_1 \\
\endmatrix
$$
where $Z \rightarrow X$ and $Z' \rightarrow X_1$ are $K$-negative extremal
divisorial
contractions, $Z \dashrightarrow Z'$ a sequence of log flips with respect to an
appropriate
log pair.

\vskip.1in

Links of type (III) (Inverses of links of type (I))

$$
\matrix
X & \dashrightarrow & Z' & & \\
& & & \searrow & \\
\downarrow & & & & X_1 \\
& & & & \downarrow \\
S & & \rightarrow & & S_1 \\
\endmatrix
$$
where $Z' \rightarrow X_1$ is a $K$-negative extremal divisorial contraction,
$X
\dashrightarrow Z'$ a sequence of log flips with respect to an appropriate log
pair and
$\rho(S/S_1) = 1$.

\vskip.1in

Links of type (IV)

$$
\matrix
X & & \dashrightarrow & & X_1 \\
\downarrow & & & & \downarrow \\
S & & & & S_1 \\
& \searrow & & \swarrow & \\
& & T & & \\
\endmatrix
$$
where $X \dashrightarrow X'$ is a sequence of log flips with respect to an
appropriate log
pair and $\rho(S/T) = \rho(S_1/T) = 1$.
\endproclaim

\proclaim{Remark 1.3}\endproclaim  All the intermediate Mori fiber spaces
$$\phi_k:X_k \rightarrow S_k$$
in the process of untwisting the given birational map by the Sarkisov program
are
the end results of
$K$-MMP over $Spec\ k$ starting from one appropriate nonsingular projective
variety
$W$.  In order to see this, we just have to take $W$ to be a common resolution
$$p_k:W \rightarrow X_k.$$
Then each $p_k$ is a process of $K$-MMP over $X_k$, and
thus a process of $K$-MMP over $Spec\ k$.

This fact that all the Mori fiber spaces in the process of the Sarkisov program
are Sarkisov related (See \S 3 for the precise definition of the Sarkisov
relation.) is automatic in the case of the genuine Sarkisov program and
implicit in the statement.  But it is the key point of understanding the log
Sarkisov program.

\vskip.1in

The strategy to untwist $\Phi$ into a composite of links is to set up a good
invariant, the Sarkisov degree $(\mu,\lambda,e)$ so that it strictly decreases
after untwisting the birational map.  That is to say, we would like to
construct
a sequence of links as below

\newpage

\vskip2in

$$
\CD
X = X_0@.\dashrightarrow @.X_1 @. \dasharrow @.X_2 @.@. @.\ \ \ \ \
\ \ \ \ \ \ \ \ \ \cdot\cdot\cdot
\ \ \ \ \ \ \ \ \ \ \ \ \ \dasharrow @.X'\\
@V\phi VV @.@V\phi_1VV @. @V\phi_2VV @.@.\cdot\cdot\cdot @. @V\phi'VV\\
S = S_0@.@. S_1 @.@.S_2 @.@.@.\cdot\cdot\cdot @.S' \\
\endCD
$$

so that the Sarkisov degree strictly decreases each time we untwist the
birational map
$$(\mu,\lambda,e) = (\mu_0,\lambda_0,e_0) > (\mu_1,\lambda_1,e_1) >
\cdot\cdot\cdot$$ and that it cannot decrease infinitely many times.

\proclaim{Definition 1.4 (Sarkisov degree)} The Sarkisov degree of a
birational map between two Mori fiber spaces
$$
\CD
X @.\overset\Phi\to\dashrightarrow @.X' \\
@V\phi VV @.@VV\phi' V \\
S @.@. S' \\
\endCD
$$
with reference to the fixed Mori fiber space $\phi':X' \rightarrow S'$ is the
triplet
$$(\mu,\lambda,e)$$ of the numbers defined below, endowed with the
lexicographical order.
\endproclaim

\vskip.1in

First we take a very ample divisor $A'$ on $S'$ and a sufficiently
divisible $\mu' \in {\Bbb N}$ such that
$${\Cal H}_{X'} = -\mu'K_{X'} + {\phi'}^*A'$$
is very ample on $X'$.  ${\Cal H}_X$ is the strict transform of ${\Cal
H}_{X'}$.

\vskip.1in

$\boxed{\mu:\text{the\ quasi-effective\ threshold}}$

\vskip.1in

The quasi-effective threshold $\mu$ is defined to be a positive rational number
s.t.
$$\mu K_X + {\Cal H}_X \equiv 0 \text{\ over\ }S.$$

Note that $\mu'$ is the quasi-effective threshold for the special case $\Phi$
being the identity map of the Mori fiber space $\phi':X' \rightarrow S'$.

In dimension 2, it is easy to see $\mu \in \frac{1}{3!}{\Bbb N}$.

\vskip.1in

$\boxed{\lambda:\text{the\ maximal\ multiplicity\ of\ an\ extremal\ ray}}$

\vskip.1in

We take a common resolution
$$
\CD
@.@.W @.@.\\
@.p\swarrow@.@.\searrow q\\
X @.@.\overset\Phi\to\dashrightarrow @.@.X' \\
@V\phi VV @.@.@.@VV\phi' V \\
S @.@.@.@. S'. \\
\endCD
$$
such that the exceptional locus of $p$ is an S.N.C. divisor $\cup E_k$.  Then
taking a general member of ${\Cal H}_{X'}$ and its strict transform ${\Cal
H}_X$ (We denote them by the same symbols ${\Cal H}_{X'}$ and ${\Cal H}_{X}$ by
abuse of notation.) and writing
$$\align
K_W &= p^*K_X + \Sigma a_kE_k\\
q^*{\Cal H}_{X'} &= p^*{\Cal H}_X - \Sigma b_kE_k,\\
\endalign$$
we define
$$\lambda := \text{max}\{\frac{b_k}{a_k}\}.$$

We remark that $\frac{1}{\lambda}$ has a more intrinsic description and is
called the
canonical threshold of $X$ with respect to ${\Cal H}_X$, i.e.,
$$\frac{1}{\lambda} = \text{max}\{c \in {\Bbb Q}_{>0};K_X + c{\Cal H}_X
\text{\ canonical}\},$$
where $K_X + c{\Cal H}_X$ being canonical means by
definition that for some (and thus for any) common resolution
$$\align
p&:W \rightarrow X\\
q&:W \rightarrow X'\\
\endalign$$
such that the exceptional locus
$\cup E_k$ of $p$ is an S.N.C. divisor $\cup E_k$ (Note that the strict
transform
${\Cal H}_W$ which is nothing but the total transform $q^*{\Cal H}_{X'}$ may
assumed to be nonsingular and cross $\cup E_k$ normally.), we have $$K_W +
c{\Cal
H}_W = p^*(K_X + c{\Cal H}_X) + \Sigma r_kE_k$$ with $$r_k \geq 0 \text{\ for\
}\forall k.$$
Thus $\lambda$ is independent of the common resolution that we
take and well-defined.  We note that when $X$ is ${\Bbb Q}$-factorial and thus
$p$ has purely one-codimensional exceptional locus, the assumption of the
exceptional locus of $p$ being an S.N.C. divisor is unnecessary.

Note also that since ${\Cal H}_X$ has no base component of codimension one,
even
when $c \geq 1$ we can regard the pair $(X,c{\Cal H}_X) = (X,\Sigma c_qB_q)$ as
a
canonical log pair with only klt singularities in the usual sense (cf. Koll\'ar
et al[16]) by taking general members $B_q \in {\Cal H}_q$ and a suitable set of
positive rational numbers
$0 < c_q < 1$ with $\Sigma c_q = c$ and that thus $K + c{\Cal H}$-MMP works as
well as
$K + \Sigma c_qB_q$-MMP.

In dimension 2, $\lambda$ is nothing but the maximal multiplicity of a general
member of the linear system ${\Cal H}_X$ (We note that the linear system ${\Cal
H}_X$ consists of the strict transforms of the complete linear system ${\Cal
H}_{X'}$ and that it may not be complete itself.) at the base points $Bs({\Cal
H}_X)$.  When $Bs({\Cal H}_X) = \emptyset$, $\lambda = 0$ by definition.

\vskip.1in

$\boxed{e:\text{the\ number\ of\ }K + \frac{1}{\lambda}{\Cal H}-\text{crepant\
divisors}}$

\vskip.1in

$$e := \left\{\aligned
&\# \{E;a(E,\frac{1}{\lambda}{\Cal H}_X) = 0\} \text{\ if\ }\lambda > 0\\
&0 \text{\ if\ }\lambda = 0
\endaligned
\right.$$

In dimension 2, $e$ is nothing but the number of the base points of the linear
system ${\Cal H}_X$ with the maximal multiplicity $\lambda$ when $Bs({\Cal
H}_X)
\neq \emptyset$, and $e = 0$ by definition when $Bs({\Cal H}_X) = \emptyset$.

\vskip.1in

Once the Sarkisov degree, which should measure the extent of untwisting, is set
up, the only other ingredient we need is a criterion to judge if the untwisting
is
completed:

\proclaim{Proposition 1.5 (N\"other-Fano criterion)}  A birational map $\Phi$
between two Mori fiber spaces is an isomorphism of Mori fiber spaces, i.e.,

$$
\CD
X @.\overset\sim\to\rightarrow @.X' \\
@V\phi VV @.@VV\phi' V \\
S @.\overset\sim\to\rightarrow @. S' \\
\endCD
$$
if $\lambda \leq \mu$ and $K_X + \frac{1}{\mu}{\Cal H}_X$ is nef.
\endproclaim

The proposition can be proved as an easy application of the Negativity Lemma or
Hodge Index Theorem and we refer the reader to Corti[4],Theorem4.2 for a proof.

\vskip.2in

$$\bold{Flowchart\ for\ Sarkisov\ Program.}$$

\vskip.1in

In the following, we present a flowchart to untwist a birational map
$$
\CD
X @.\overset\Phi\to\dashrightarrow @.X' \\
@V\phi VV @.@VV\phi' V \\
S @.@. S' \\
\endCD
$$
between two Mori fiber spaces.

\vskip.1in

We START.

The first question to ask is:
$$\lambda > \mu ?$$
According to whether the answer to this question is YES or NO, we proceed
separately into the case $\lambda > \mu$ or into the case $\lambda \leq \mu$.

$$\boxed{\text{Case}:\lambda \leq \mu}$$

If $\lambda \leq \mu$, then the next question to ask is:
$$K_X + \frac{1}{\mu}{\Cal H}_X \text{\ nef}?$$
If the answer to this question is YES, then $K_X + \frac{1}{\mu}{\Cal H}_X$ is
nef
and $\lambda \leq \mu$ by the case assumption.  Thus the N\"other-Fano
criterion applies to conclude $\Phi$ is an isomorphism of Mori fiber spaces.
This leads to an
$$\text{END.}$$
If $K_X + \frac{1}{\mu}{\Cal H}_X$ is not nef, then we construct as follows a
normal projective variety $T$ dominated by $S \rightarrow T$ s.t. $K_X +
\frac{1}{\mu}{\Cal H}_X$ is not relatively nef over $T$ and $\rho(X/T) = 2$, so
that we run $K + \frac{1}{\mu}{\Cal H}$-MMP over $T$ to have an untwisting
link.

We pick a $K_X + \frac{1}{\mu}{\Cal H}_X$-negative extremal ray $P$ of
$\overline{NE}(X/Spec\ k)$ s.t. the span
$F := P + R$ is a 2-dimensional extremal face, where $R$ is the $K_X$-negative
extremal ray giving the Mori fiber space $\phi:X \rightarrow S$.  $F$ is $K_X +
(\frac{1}{\mu} - \epsilon){\Cal H}_X$-negative for $0 < \epsilon << 1$, thus we
have the contraction morphism $cont_F:X \rightarrow T$ to obtain $T$.  Since $F
\supset R$, $cont_F$ factors through $S$ and by construction $T$ satisfies all
the desired conditions.

Now we
$$\text{Run\ }K + \frac{1}{\mu}{\Cal H}-\text{MMP\ over\ }T.$$

We reach either a minimal model or a Mori fiber space (with respect to $K +
\frac{1}{\mu}{\Cal H}$ and over $T$).

First we show that it is
$$\text{IMPOSSIBLE\ to\ reach\ a\ minimal\ model!}$$
Suppose we did.  Then according to whether the first nonflipping contraction is
divisorial or not, we should have two different diagrams as follows:

$$
\matrix
X & \dashrightarrow & Z' & & \\
\downarrow & & & \searrow & \\
S & & & & X_1 \\
& \searrow & & & \downarrow \\
& & T & \overset\sim\to\rightarrow & S_1 \\
\endmatrix
$$

$$
\matrix
X & & \dashrightarrow & & Z' = X_1 \\
\downarrow & & & & \downarrow \\
S & & & & S_1 \\
& \searrow & & \swarrow & \\
& & T & & \\
\endmatrix
$$
We take a general curve $\Sigma_1 \in Hilb(X_1/T)$ away from the locus of
indeterminacy of the birational map $X_1 \dashrightarrow X$ (i.e., in the first
case the union of the image of the exceptional divisor of the divisorial
contraction and all the flipped curves and in the second case all the flipped
curves).  $\Sigma_1$ can be considered to lie on $X$ and since $\Sigma_1$ is
general we conclude $\phi(\Sigma_1)$ is a curve (not a point), which implies
$(K_X +
\frac{1}{\mu}{\Cal H}_X) \cdot \Sigma_1 < 0$.  But then
$$0 > (K_X + \frac{1}{\mu}{\Cal H}_X) \cdot \Sigma = (K_{X_1} +
\frac{1}{\mu}{\Cal
H}_{X_1}) \cdot \Sigma_1 \geq 0,$$
a contradiction!

Next suppose we
$$\text{Reach\ a\ Mori\ fiber\ space\ }X_1 \rightarrow S_1.$$

Then the next question to ask just in order to separate the types of links is:
$$\text{Is\ the\ first\ nonflipping\ contraction\ divisorial\ ?}$$

If the answer is YES, the $K + \frac{1}{\mu}{\Cal H}$-MMP consists of a
sequence
of $K + \frac{1}{\mu}{\Cal H}$-flips $X \dashrightarrow Z'$ followed by a $K +
\frac{1}{\mu}{\Cal H}$-negative divisorial contraction $Z' \rightarrow X_1$.
Since $\rho(X_1/T) = 1$, $\phi_1:X_1 \rightarrow S_1 = T$ is a $K_{X_1} +
\frac{1}{\mu}{\Cal H}_{X_1}$-negative and thus $K_{X_1}$-negative fiber space.

$$
\matrix
X & \dashrightarrow & Z' & & \\
\downarrow & & & \searrow & \\
S & & & & X_1 \\
& \searrow & & & \downarrow \\
& & T & \overset\sim\to\rightarrow & S_1 \\
\endmatrix
$$

If the answer is NO, then the $K + \frac{1}{\mu}{\Cal H}$-MMP consists of a
sequence of $K + \frac{1}{\mu}{\Cal H}$-flips $X \dashrightarrow Z'$ followed
by a
$K_{X_1} + \frac{1}{\mu}{\Cal H}_{X_1}$-negative and thus $K_{X_1}$-negative
fibering contraction $\phi_1:Z' = X_1 \rightarrow S_1$.  Since $\rho(X_1/T) =
\rho(X/T) = 2$, we have $\rho(S_1/T) = 1$.

$$
\matrix
X & & \dashrightarrow & & Z' = X_1 \\
\downarrow & & & & \downarrow \\
S & & & & S_1 \\
& \searrow & & \swarrow & \\
& & T & & \\
\endmatrix
$$

We claim in both cases $X_1$ has only terminal singularities.  (${\Bbb
Q}$-factoriality of $X_1$ is automatic from construction.)  Let $I$ be the
locus
of indeterminacy of the birational map $X_1 \dashrightarrow X$.  If $E$ is a
discrete valuation whose center on $X_1$ is not contained in $I$ (and has
codimension $\geq 2$), then
$$a(E,X_1,\emptyset) = a(E,X,\emptyset) > 0.$$
If the center of $E$ on $X_1$ is contained in $I$, then
$$\align
a(E,X_1,\emptyset) &\geq a(E,X_1,\frac{1}{\mu}{\Cal H}_{X_1})\\
&> a(E,X,\frac{1}{\mu}{\Cal H}_X) \geq 0.\\
\endalign$$
Thus we have the claim.

Therefore, we have a link of type (III) in the former case and a link of type
(IV)
in the latter.

Moreover, since $K_{X_1} + \frac{1}{\mu}{\Cal H}_{X_1}$ is negative over $S_1$,
we
conclude in both cases
$$\mu_1 < \mu.$$

Therefore, after untwisting $\Phi$ by a link of type (III) or type (IV), we go
back to the START with strictly decreased quasi-effective threshold.

\vskip.2in

$$\boxed{\text{Case}:\lambda > \mu}$$

\vskip.1in

In this case we
$$\text{Take\ a\ maximal\ divisorial\ blow\ up\ }p:Z \rightarrow X,$$
with respect to $K_X + \frac{1}{\lambda}{\Cal H}_X$, i.e., $p$ is a projective
morphism from $Z$ with only ${\Bbb Q}$-factorial terminal singularities
s.t.

i) $\rho(Z/X) = 1$,

ii) the exceptional locus of $p$ is a prime divisor $E$, and

iii) $p$ is $K + \frac{1}{\lambda}{\Cal H}$-crepant, i.e.,
$$K_Z + \frac{1}{\lambda}{\Cal H}_Z = p^*(K_X + \frac{1}{\lambda}{\Cal H}_X).$$

\proclaim{Proposition 1.5} A maximal divisorial blow up $p:Z \rightarrow X$
with respect to $K_X + \frac{1}{\lambda}{\Cal H}_X$ exists.
\endproclaim

We remark that the exceptional divisor $E$ of $p$ is necessarily one of the $K
+
\frac{1}{\lambda}{\Cal H}$-crepant divisors $\{E_1, E_2,
\cdot\cdot\cdot, E_e\}$ counted for the number $e$.  As long
as we require $Z$ to have only terminal singularities, we can't quite
specify which $E_i$ would be the exceptional divisor.  On the other hand,
if we allow $Z$ to have canonical singularities, for each $E_i$ we can
construct a maximal blow up $p_i:Z_i \rightarrow X$ (allowing $Z_i$ to
have canonical singularities) with the exceptional divisor being $E_i$.

\demo{Proof}\enddemo

Take a resolution $Y \rightarrow X$ s.t.

a) the exceptional locus is a divisor with only S.N.C.,

b) $Y$ dominates $X'$ so that the strict transform ${\Cal H}_Y$ coincides
with the total transform of ${\Cal H}_{X'}$ and that a general
member ${\Cal H}_Y$ is smooth and crosses normally with the exceptional
locus.

We run the $K + \frac{1}{\lambda}{\Cal H}$-MMP over $X$ to get a minimal
model $f:(Z',\frac{1}{\lambda}{\Cal H}_{Z'}) \rightarrow
(X,\frac{1}{\lambda}{\Cal H}_X)$.  As before, it is easy to see that $Z'$
has only ${\Bbb Q}$-factorial terminal singularities.  Since both $Z'$
and $X$ are ${\Bbb Q}$-factorial, the exceptional locus of $f$ is purely
one-codimensional.  An easy application of the Negativity Lemma
(Shokurov[30]Koll\'ar et al[16]Corti[4]) shows that the exceptional
locus is actually $\cup_{i=1}^e E_i$ and that $f$ is $K +
\frac{1}{\lambda}{\Cal H}$-crepant, i.e., $K_{Z'} +
\frac{1}{\lambda}{\Cal H}_{Z'} = f^*(K_X + \frac{1}{\lambda}{\Cal H}_X)$.

Now we run the $K$-MMP starting from $Z'$ over $X$ ending necessarily
with a divisorial contraction $p:Z \rightarrow X$.  It is immediate that
$p:Z \rightarrow X$ is a maximal divisorial blow up with respect to $K_X
+ \frac{1}{\lambda}{\Cal H}_X$.  (If we want to specify the exceptional
divisor $E_i$ allowing $Z$ to have canonical singularities, then we run
$K + \frac{1}{\lambda}{\Cal H} + \epsilon \Sigma_{j\neq i}E_j$-MMP ($0 <
\epsilon << 1$) instead.)

We also remark that in order to construct just one maximal divisorial blowup
we can start from any common resolution $Y$ which may not satisfy a) as long
as $X$ is ${\Bbb Q}$-factorial and thus the exceptional locus of $Y
\rightarrow X$ is purely one codimensional.

\vskip.1in

There is another method called the ``Nef Threshold Method" to construct a
maximal divisorial blow up by M. Reid.

We construct a chain of

$Y_i$: 3-folds with only ${\Bbb Q}$-factorial terminal singularities
projective over $X$ $(Y_0 = Y)$,

${\Cal H}_{Y_i}$: the strict transforms of ${\Cal H}_Y$,

$\lambda_i$: (a non-decreasing sequence of) nonnegative rational numbers

s.t.

a) $\lambda_iK_{Y_i} + {\Cal H}_{Y_i}$ is a supporting function of a face
containing a $K_{Y_i}$-negative extremal ray $R_i$ of $\overline{NE}(Y_i/X)$,
i.e., $\lambda_iK_{Y_i} + {\Cal H}_{Y_i}$ is relatively nef over $X$ and
$$(\lambda_iK_{Y_i} + {\Cal H}_{Y_i})^{\perp} \cap \overline{NE}(Y_i/X)
\supset R_i,$$

b) either $Y_i \rightarrow Y_{i+1}$ is a divisorial contraction of $R_i$ or
the flip $Y_i \dashrightarrow Y_{i+1}$, and

c) the chain ends with a divisorial contraction $p:Z = Y_n \rightarrow X$ of
an extremal ray $R_n$
$$(\lambda_iK_{Y_i} + {\Cal H}_{Y_i})^{\perp} \cap \overline{NE}(Y_i/X)
\supset R_n.$$

Then it is easy to see that $p:Z \rightarrow X$ is a maximal divisorial blow
up with respect to $K_X + \frac{1}{\lambda}{\Cal H}_X$ and $\lambda_n =
\lambda$.

We construct inductively.

Suppose we have succeeded constructing the chain up to the $i$-th stage.
Consider the nef threshold $\lambda_i$ of ${\Cal H}_{Y_i}$ with respect to
$K_{Y_i}$
$$\lambda_i := \text{sup}\{\nu;\nu K_{Y_i} + {\Cal H}_{Y_i} \text{\ relatively\
nef\ over\ }X\}.$$
Remark that since $\lambda_{i-1}K_{Y_{i-1}} + {\Cal H}_{Y_{i-1}}$ is
relatively nef over $X$ and the contraction of $R_{i-1}$ is
$\lambda_{i-1}K_{Y_{i-1}} + {\Cal
H}_{Y_{i-1}}$-trivial, $\lambda_{i-1}K_{Y_i} + {\Cal H}_{Y_i}$ is also
relatively nef over $X$ and thus $\lambda_i \geq \lambda_{i-1}$.

We claim that $\lambda_i$ is rational and that there exists a
$K_{Y_i}$-negative extremal ray $R_i$ s.t.
$$(\lambda_iK_{Y_i} + {\Cal H}_{Y_i})^{\perp} \cap \overline{NE}(Y_i/X) \supset
R_i.$$
Instead of applying the Rationality Theorem (KaMaMa[11],Theorem4-1-1) whose
proof only applies to ample divisors, we use a result of Kawamata[8] on the
boundedness of lengths of the extremal rational curves to the relatively nef
divisor $\lambda_{i-1}K_{Y_i} + {\Cal H}_{Y_i}$.  (This idea was communicated
to
us by J. McKernan.  See KeMaMc[12].)  First from the definition and
the Cone Theorem, we have
$$\align
\lambda_i &= \lambda_{i-1} + \text{inf}\{\frac{(\lambda_{i-1}K_{Y_i} + {\Cal
H}_{Y_i})\cdot l}{-K_{Y_i} \cdot l}\};l:K_{Y_i}-\text{negative\ extremal\
rays}\}\\
&= \text{inf}\{\frac{{\Cal H}_{Y_i}\cdot l}{-K_{Y_i} \cdot
l}\};l:K_{Y_i}-\text{negative\ extremal\ rays}\}.\\
\endalign$$
The result of Kawamata tells us that for each $K_{Y_i}$-negative extremal ray
$l$, there exists a rational curve $L_l$ which generates $l = {\Bbb R}_+[L_l]$
and $0 < -K_{Y_i} \cdot L_l \leq 2 \cdot \dim X$.  Thus if $q_i$ is the
${\Bbb Q}$-factorial index of $Y_i$
$$q_i := \text{min}\{z \in {\Bbb N};zD \text{\ is\ Cartier\ for\ all\ integral\
Weil\ divisors\ }D \text{\ on\ }Y_i\},$$
(which coincides with the index $r_i$ of $K_{Y_i}$ in dimension 3), then
$$\frac{{\Cal H}_{Y_i}\cdot l}{-K_{Y_i} \cdot
l} \in \frac{1}{(r_i \cdot 2 \text{dim}X)!q_i}{\Bbb Z}_{\geq 0}.$$
Therefore, ``inf" is actually attained as the minimum for some
$K_{Y_i}$-negative extremal ray $R_i$ and for this $R_i$ we have
$$R_i \subset (\lambda_iK_{Y_i} + {\Cal H}_{Y_i})^{\perp} \cap
\overline{NE}(Y_i/X).$$
As for $Y_{i+1}$, we take either the divisorial contraction $Y_i \rightarrow
Y_{i+1}$ of $R_i$ or the flip $Y_i \dashrightarrow Y_{i+1}$ of $R_i$.

\vskip.2in

We go back to the discussion of the flowchart.

\vskip.1in

Now we
$$\text{Run\ }K + \frac{1}{\lambda}{\Cal H}-\text{MMP\ over\ }S.$$
A priori we reach either a minimal model or a Mori fiber space (with respect
to $K + \frac{1}{\lambda}{\Cal H}$ over $S$).

First we show that it is
$$\text{IMPOSSIBLE\ to\ reach\ a\ minimal\ model!}$$

Suppose we did.

Then according to whether the first nonflipping contraction is divisorial or
not we should have two different diagrams:

$$
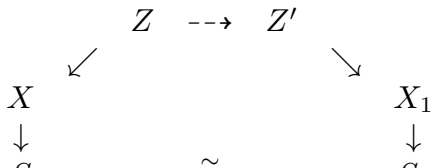
\matrix
& & Z & \dashrightarrow & Z' & & \\
& \swarrow & & & & \searrow & \\
X & & & & & & X_1 \\
\downarrow & & & & & & \downarrow \\
S & & & \overset\sim\to\leftarrow & & & S_1 \\
\endmatrix
$$

$$
\matrix
& & Z & \dashrightarrow & X_1 \\
& \swarrow & & & \\
X & & & & \downarrow \\
\downarrow & & & & \\
S & & \leftarrow & & S_1 \\
\endmatrix
$$

We take a general curve $\Sigma \in Hilb(X/S)$ away from $p(E)$ (and thus can
be considered to lie on $Z$) and away from all the flipping curves (and thus
can be considered to lie on $Z'$).

In the first case, we have
$$\align
0 &\leq (K_{X_1} + \frac{1}{\lambda}{\Cal H}_{X_1}) \cdot q_*\Sigma\\
&= \{(K_{Z'} + \frac{1}{\lambda}{\Cal H}_{Z'}) - a E_q\} \cdot \Sigma\ (a >
0)\\
&\leq (K_{X} + \frac{1}{\lambda}{\Cal H}_{X}) \cdot \Sigma \\
&< (K_{X} + \frac{1}{\mu}{\Cal H}_{X}) \cdot \Sigma = 0,\\
\endalign$$
a contradiction!

In the second case, we have
$$\align
0 &\leq (K_{Z'} + \frac{1}{\lambda}{\Cal H}_{Z'}) \cdot \Sigma \\
&= (K_{X} + \frac{1}{\lambda}{\Cal H}_{X}) \cdot \Sigma \\
&< (K_{X} + \frac{1}{\mu}{\Cal H}_{X}) \cdot \Sigma = 0,\\
\endalign$$
again a contradiction!

\vskip.1in

Next suppose we
$$\text{Reach\ a\ Mori\ fiber\ space\ }X_1 \rightarrow S_1.$$

Then the next question to ask just in order to separate the types of links is:
$$\text{Is\ the\ first\ nonflipping\ contraction\ divisorial\ ?}$$

If the answer is YES, the $K + \frac{1}{\lambda}{\Cal H}$-MMP consists of a
sequence of $K + \frac{1}{\lambda}{\Cal H}$-flips $X \dashrightarrow Z'$
followed by a $K + \frac{1}{\lambda}{\Cal H}$-negative contraction $Z'
\rightarrow X_1$.  Since $\rho(X_1/S) = 1$, $\phi_1:X_1 \rightarrow S_1 = S$
is a $K_{X_1} + \frac{1}{\lambda}{\Cal H}_{X_1}$-negative and thus
$K_{X_1}$-negative fiber space.

$$
\matrix
& & Z & \dashrightarrow & Z' & & \\
& \swarrow & & & & \searrow & \\
X & & & & & & X_1 \\
\downarrow & & & & & & \downarrow \\
S & & & \overset\sim\to\leftarrow & & & S_1 \\
\endmatrix
$$

We note that the exceptional divisors $E$ and $E_q$ are distinct, since
 otherwise $X$ and $X_1$ are isomorphic in codimension one, which would imply
$X$ and $X_1$ are indeed isomorphic over $S$, but then while $E_q$ is NOT
$K_{X_1}
 + \frac{1}{\lambda}{\Cal H}_{X_1}$-crepant
 $E$ is $K_X + \frac{1}{\lambda}{\Cal H}_X$-crepant, absurd!

\vskip.1in

If the answer is NO, then the $K + \frac{1}{\lambda}{\Cal H}$-MMP consists of
a sequence of $K + \frac{1}{\lambda}{\Cal H}$-flips $X \dashrightarrow Z'$
followed by a $K_{X_1} + \frac{1}{\lambda}{\Cal H}_{X_1}$-negative and thus
$K_{X_1}$-negative fibering contraction $Z' = X_1 \rightarrow S_1$.  Since
$\rho(X_1/S) = \rho(Z/S) = 2$, we have $\rho(S_1/T) = 1$.

$$
\matrix
& & Z & \dashrightarrow & X_1 \\
& \swarrow & & & \\
X & & & & \downarrow \\
\downarrow & & & & \\
S & & \leftarrow & & S_1 \\
\endmatrix
$$

In both cases, $X_1$ has only ${\Bbb Q}$-factorial terminal singularities and
thus we have a link of type (II) or a link of type (I), respectively.

\vskip.1in

Now we study how the Sarkisov degree $(\mu,\lambda,e)$ changes after
untwisting by a link of type (II) or type (I).

We claim that
$$\mu_1 \leq \mu$$
with equality holding only if
$$\align
\text{either\ }&\dim S_1 > \dim S\\
\text{or\ }&\dim S_1 = \dim S\text{\ and\ }\psi_1\text{\ is\
square\ },\\
\endalign$$
i.e.,
$$
\CD
X @.\overset\psi_1\to\dashrightarrow @.X_1 \\
@V\phi VV @.@VV\phi_1 V \\
S @.\overset\pi\to\leftarrow @. S' \\
\endCD
$$
$\pi$ is a birational morphism and $\psi_{\eta}:X_{\eta} \dashrightarrow
(X_1)_{\eta}$ is an isomorphism, where $\eta$ is the generic point of $S$.

First by definition of $\lambda$ and the assumption of this case $\lambda >
\mu$, it follows that
$$p^*(K_X + \frac{1}{\mu}{\Cal H}_X) = K_Z + \frac{1}{\mu}{\Cal H}_Z + bE$$
for some $b > 0, b \in {\Bbb Q}$.

We take a general curve $\Sigma_1 \in Hilb(X_1/S_1)$ away from the locus of
indeterminacy of the birational map $X_1 \dashrightarrow Z$ (i.e., in the
case of a link of type (II) the union of $q(E_q)$ and all the flipped
curves and in the case of a link of type (I) the union of all the flipped
curves).  Then $\Sigma_1$ can be considered to lie on $Z$ and
$$\align
0 &= (K_X + \frac{1}{\mu}{\Cal H}_X) \cdot p_*\Sigma_1\\
&= (K_Z + \frac{1}{\mu}{\Cal H}_Z + bE) \cdot \Sigma_1\\
&\geq (K_Z + \frac{1}{\mu}{\Cal H}_Z) \cdot \Sigma_1\\
&= (K_{X_1} + \frac{1}{\mu}{\Cal H}_{X_1}) \cdot \Sigma_1,\\
\endalign$$
which implies
$$\mu_1 \leq \mu.$$
Moreover, if $\mu_1 = \mu$ and $\dim S = \dim S_1$ (which implies
that $\pi:S_1 \rightarrow S$ is a birational morphism, since both field
extensions $k(X)/k(S)$ and $k(X_1) = k(X)/k(S_1)$ are algebraically closed),
then $E \cdot \Sigma_1 = 0$, which is equivalent to saying
\linebreak
$\phi_1(\text{the\ strict\ transform\ of\ }E) \neq S_1$.  Therefore,
$\psi_1$ is square.

\vskip.1in

We also claim that
$$\lambda_1 \leq \lambda$$
and
$$\text{if\ }\lambda_1 = \lambda\text{\ then\ }e_1 < e.$$

First $(X_1,\frac{1}{\lambda}{\Cal H}_{X_1})$ is canonical, since it is
obtained from a canonical pair $(Z,\frac{1}{\lambda}{\Cal H}_Z)$ through $K +
\frac{1}{\lambda}{\Cal H}$-MMP.  Thus $\lambda_1 \leq \lambda$.  (Note that in
general canonicality may not be preserved when we contract a component of the
boundary $B$ through $K + B$-MMP.  But in our case, ${\Cal H}$'s are the
strict transforms of one unique base point free system and thus canonicality
is preserved.)

Moreover, if $\lambda_1 = \lambda$, then in the case of untwisting by a link
of type (II)
$$K_{Z'} + \frac{1}{\lambda}{\Cal H}_{Z'} = q^*(K_{X_1} +
\frac{1}{\lambda}{\Cal H}_{X_1}) + aE_q (a > 0)$$
implies $E_q$ is not a $K_{X_1} + \frac{1}{\lambda}{\Cal H}_{X_1}$-crepant
divisor (and $E$ is a divisor on $X_1$ and thus not exceptional) and thus
$$e_1 \leq e - 1 < e.$$
In the case of untwisting
 by a link of type (I) $E$ is a divisor on $X_1$ (and thus not exceptional)
and thus we have the same conclusion.

\vskip.1in

Therefore, after untwisting by a link of type (II) or type (I), we go back to
the START with strictly decreased Sarkisov degree.

\vskip.1in

The ``visualization" of the flowchart can be found at the end of the paper.

\newpage

\subhead{\S 2. Termination of Flowchart}\endsubhead

\vskip.1in

In this section, we discuss the termination of the flowchart for the Sarkisov
program, i.e., the problem of showing that there is no infinite loop in the
flowchart and thus after a finite number of untwisting it gives a factorization
of any given birational map between two Mori fiber spaces.  Once we have
(Log)-MMP in dimension $n$ the key points of showing termination for
Sarkisov program for $n$-folds are:

i) Discreteness of the quasi-effective thresholds $\mu$, which follows from the
boundedness of ${\Bbb Q}$-Fano $d$-folds $d \leq n$, and

ii) Corti[4]'s ingeneous argument to reduce the probelm to
$S_n(\text{Local})$ when the quasi-effective threshold stabilizes.

In dimension 3, where we have all the necessary ingredients, the termination of
the flowchart is a theorem by Corti[4].  The argument here is a modification of
Corti[4] following a slightly simplified flowchart in the previous section.
We restrict ourselves to dimension 3 in the following presentation, but we
carry the argument so that it works almost verbatim in arbitrary dimension
(once all the necessary but still conjectural ingredients are established).

\vskip.1in

\proclaim{Claim 2.1} There is no infinite number of untwisting (successive or
unsuccessive) by the links under the case $\lambda \leq \mu$.
\endproclaim
\demo{Proof}\enddemo Suppose there are infinitely many links (successive or
unsuccessive)
$$
\CD
X_i @.\overset\psi_i\to\dashrightarrow @.X_{i+1} \\
@V\phi_i VV @.@VV\phi_{i+1} V \\
S_i @.@. S_{i+1} \\
\endCD
$$
under the case $\lambda \leq \mu$.  Note that in the case $\lambda \leq \mu$
we have $\dim S_i \geq 1$ (unless $\Phi_i$ becomes an isomorphism of
Mori fiber spaces).  When
$\dim S_i = 2$,
$l$ being a rational curve which is a general fiber of $\phi_i$, we have
$$\align
K_{X_i} \cdot l &= -2\\
(\mu K_{X_i} + {\Cal H}_{X_i}) \cdot l &= 0,\\
\endalign$$
which implies
$$\mu \in \frac{1}{2}{\Bbb N}.$$
When $\dim S_i = 1$, we can take a rational curve $l$ in a general fiber
which is a Del Pezzo surface s.t.
$$\align
K_{X_i} \cdot l &= -1, -2 \text{\ or\ }-3\\
(\mu K_{X_i} + {\Cal H}_{X_i}) \cdot l &= 0,\\
\endalign$$
which implies
$$\mu \in \frac{1}{3!}{\Bbb N}.$$
Since after any link in the case $\lambda \leq \mu$ the quasi-effective
threshold strictly decreases and it does not increase afyter any link in any
case, we then have a strictly decreasing sequence in $\frac{1}{3!}{\Bbb N}$
$$\mu \geq \mu_1 > \mu_2 \cdot\cdot\cdot > 0,$$
a contradiction!

In general, we only have to use the boundedness of ${\Bbb Q}$-Fano $d$-folds
for $d \leq n - 1$ to derive the discreteness of $\mu$ and thus a
contradiction to establish this claim.

\vskip.1in

\proclaim{Claim 2.2} There is no infinite (successive) sequence of untwisting
by the links under the case $\lambda > \mu$ with stationary quasi-effective
threshold.
\endproclaim

\demo{Proof}\enddemo This is the heart of the ingeneous argument by
Corti[4].  Suppose there is such an infinite sequence
$$
\CD
X = X_0@.\dashrightarrow @.X_1 @. \dasharrow @.X_2 @.\dashrightarrow @. @.\ \ \
\
\
\ \ \ \ \ \ \ \ \ \cdot\cdot\cdot
\ \ \ \ \ \ \ \ \ \ \ \ \ \dasharrow @.X_{k} @.\overset\psi_k\to\dashrightarrow
@.X_{k+1} @. \cdot\cdot\cdot \\
@V\phi VV @.@V\phi_1VV @. @V\phi_2VV @.@.\cdot\cdot\cdot @.
@V\phi_k VV @. @V\phi_{k+1}VV @. \cdot\cdot\cdot \\
S = S_0@.\leftarrow @. S_1 @.\leftarrow @.S_2 @.\leftarrow @.@.
\ \ \ \ \ \ \ \ \ \ \ \cdot\cdot\cdot
\ \ \ \ \ \ \ \ \ \ \ \ \ \ \ \ \leftarrow @.S_k @.\leftarrow @. S_{k+1} @.
\cdot\cdot\cdot
\\
\endCD
$$
Since $\mu_k = \mu_{k+1}$ for each $k$ by assumption, we have
$$\align
\text{either\ } &\dim S_{k+1} > \dim S_k\\
\text{or\ } &\dim S_{k+1} = \dim S_k \text{\ and\ }\psi_k\text{\ is\
square.}\\
\endalign$$
The first cannot happen infinitely many times, thus we may assume we have the
second case for all $k$.  Note that $\dim S_k \geq 1$, since if
$\dim S_k = \dim S_{k+1} = 0$ then $\psi_k$ being square would imply
$\psi_k$ is an isomorphism of Mori fiber spaces, which is absurd!

We also know that $\{\lambda
_k\}$ is a nonincreasing sequence and since if $\lambda_k =
\lambda_{k+1}$ then $e_{k+1} < e_k$, the value of $\lambda_k$ cannot be
stationary.  Therefore, we have a sequence
$$\{\frac{1}{\lambda_k}\}\ \ (\frac{1}{\lambda_k} < \frac{1}{\mu_k} =
\frac{1}{\mu_0})$$
which accumulates from below to (but never equals)
$$\alpha \leq \frac{1}{\mu_0}.$$

Step 1. We claim $(X_k,\alpha{\Cal H}_{X_k})$ and $(Z_k,\alpha{\Cal
H}_{Z_k})$ ($p_k:Z_k \rightarrow X_k$ is a maximal divisorial blowup with
respect to $K + \frac{1}{\lambda}{\Cal H}$) have only log canonical
singularities for $k$ sufficiently large (and thus we may assume this
holds for $\forall k$).

Let $\alpha_k (> \frac{1}{\lambda_k})$ be the log canonical threshold of the
pair $X_k$ with respect to ${\Cal H}_k$.  If $\alpha > \alpha_k$ for infinitely
many $k$'s, then there is a strictly increasing subsequence $\{\alpha_l\}$ of
log canonical thresholds accumulating to $\alpha$.  This contradicts
$S_3(\text{Local})$ proved by Alexeev[1] (cf.Koll\'ar et al[16]).  The same
argument applies to $(Z_k,\alpha{\Cal H}_{Z_k})$.

\vskip.1in

Every link $X_k \dashrightarrow X_{k+1}$ is an outcome of $K +
\frac{1}{\lambda_k}{\Cal H}$-MMP over $S_k$ (after taking a maximal
divisorial blowup $p_k:Z_k \rightarrow X_k$) consisting of a finite number
of $K + \frac{1}{\lambda}{\Cal H}$-flips
$$Z_k = Z_k^0 \overset t^0\to\dashrightarrow Z_k^1
\overset t^1\to\dashrightarrow Z_k^2 \cdot\cdot\cdot
\overset t^{m-1}\to\dashrightarrow Z_k^m,$$
possibly followed by a divisorial contraction $q_k^m:Z_k^m \rightarrow
X_k^{m+1} = X_{k+1}$ (otherwise $Z_k^m = X_k^{m+1}$).

\vskip.1in

Step 2. We claim that every step

$$
\matrix
Z_k^i & & \dashrightarrow & & Z_k^{k+1} \\
& \searrow & & \swarrow & \\
& & X_k^{i+1} & & \\
\endmatrix
$$

is a step of $K + \alpha{\Cal H}$-MMP.

We prove this by induction on $i$.

First note that since $\alpha > c_k$, we have
$$K_{Z_k} + \alpha{\Cal H}_{Z_k} = {p_k^0}^*(K_{Z_k} + \alpha{\Cal H}_{Z_k}) -
aE_k (a > 0).$$
Therefore, we have
$$(K_{Z_k^0}
 + \alpha{\Cal H}_{Z_k^0}) \cdot P_k^0 > 0,$$
$P_k^0$ being the extremal ray giving rise to the morphism $p_k^0$.

Suppose we have
$$(K_{Z_k^i}
 + \alpha{\Cal H}_{Z_k^i}) \cdot P_k^i > 0,$$
$P_k^i$ being the extremal ray giving rise to the morphism $p_k^i$.

Note that $K_{Z_k^i}
 + \alpha{\Cal H}_{Z_k^i}$ is never relatively nef over $S_k$.  We see this as
follows: First $\alpha \leq \frac{1}{\mu_k} = \frac{1}{\mu_1}$.  Suppose
$\alpha = \frac{1}{\mu_k}.$  Then
$$K_{Z_k^i} + \alpha{\Cal H}_{Z_k^i} \equiv_{S_k} - a(\text{the\ strict\
transform\ of\ }E_k) (a > 0)$$
is never relatively nef over $S_k$.  Suppose $\alpha < \frac{1}{\mu_k}$.  Then
by taking a general curve $\Sigma \in Hilb(X_k/S_k)$ away from the locus of
indeterminacy of the birational map $X_k \dashrightarrow Z_k^i$ (which thus can
be considered to lie on $Z_k^i$) we have
$$
\align
(K_{Z_k^i} + \alpha{\Cal H}_{Z_k^i}) \cdot \Sigma &= (K_{X_k} + \alpha{\Cal
H}_{X_k}) \cdot \Sigma\\
&< (K_{X_k}
 + \frac{1}{\mu_k}{\Cal H}_{X_k}) \cdot \Sigma = 0.\\
\endalign$$
This implies that
$$(K_{Z_k^i} + \alpha{\Cal H}_{Z_k^i}) \cdot Q_k^i < 0$$
for the other extremal ray $Q_k^i$ of 2-dimensional cone
$\overline{NE}(Z_k^i/S_k)$.  This proves the claim.

\vskip.1in

A consequence of this claim is that (cf.KaMaMa[11],Proposition 5-1-11
)
$$a(\nu,X_1,\alpha{\Cal H}_{X_1}) \leq a(\nu,X_k,\alpha{\Cal H}_{X_k})$$
for any discrete valuation $\nu$ of $k(X)$ and the strict inequality holds iff
$\psi_i$ is not an isomorphism at the center of $\nu$ on $X_i$ for some $i <
k$.

\vskip.1in

Step 3. We claim that $(X_k,\alpha{\Cal H}_{X_k})$ has purely log terminal
singularities for $k$ sufficiently large (and thus we may assume this holds
for $\forall k$).

Assume to the contrary that there exists infinitely many $k$ s.t.
$(X_k,\alpha{\Cal H}_{X_k})$ is not purely log terminal, which is equivalent
to saying by the consequence above that for all $k$ there exists a valuation
$\nu_k$ of $k(X)$ with
$$a(\nu_k,X_k,\alpha{\Cal H}_{X_k}) = -1,$$
which implies again by the consequence that
$$a(\nu_k,X_1,\alpha{\Cal H}_{X_1}) = -1$$
and that at the center $z(\nu_k,X_1)$ of $\nu_k$ on $X_1$, the birational map
$\psi_{k-1} \circ \cdot\cdot\cdot \psi_2 \circ \psi_1:X_1 \dashrightarrow
X_k$ is an isomorphism.  Thus the local (w.r.t. Zariski topology) canonical
thresholds are the same
$$c(z(\nu_k,X_k),X_k,{\Cal H}_{X_k}) = c(z(\nu_k,X_1),X_1,{\Cal H}_{X_1}).$$
On the other hand, by definition
$$\frac{1}{\lambda_k} \leq c(z(\nu_k,X_k),X_k,{\Cal H}_{X_k})$$ and since
$K_{X_k} + \alpha{\Cal H}_{X_k}$ is not canonical at the center
$z(X_k,X_k)$, we have
$$c(z(\nu_k,X_k),X_k,{\Cal H}_k) < \alpha.$$
Therefore,
$$\frac{1}{\lambda_k} \leq c(z(\nu_k,X_1),X_1,{\Cal H}_{X_1}) < \alpha.$$
But $\{\frac{1}{\lambda_k}\}$ is a nondecreasing and nonstationary sequence
converging to $\alpha$ and it is easy to see the set $\{c(x,X_1,{\Cal
H}_{X_1});x \in X\}$ is finite, a contradiction!

\vskip.1in
We remark that the valuations of $k(X)$ corresponding to the $E_k$'s are all
distinct.  In fact, suppose $E_i$ and $E_j$ coincide, and thus $Z_i$ and $Z_j$
are isomorphic in a neighborhood of $E_i$ and $E_j$, which would imply
$$a(E_i,X_i,\alpha{\Cal H}_{X_i}) = a(E_j,X_j,\alpha{\Cal H}_{X_j}).$$
On the other hand, from Step 2 we have
$$a(E_i,X_i,\alpha{\Cal H}_{X_i}) < a(E_j,X_j,\alpha{\Cal H}_{X_j}),$$
a contradiction!

\vskip.1in

Finally we conclude the proof of the claim as follows: From Step 3 we may
assume that $(X_1,\alpha
{\Cal H}_{X_1})$ has only purely log terminal singularities.  But on the other
hand, for infinitely many $E_k$ with distinct corresponding discrete valuations
$$a(E_k,X_1,\alpha{\Cal H}_{X_1}) \leq a(E_k,X_k,\alpha{\Cal H}_{X_k}) < 0,$$
a contradiction!

In general, we only need $S_n(\text{Local})$ to carry out the argument for
this claim.

\vskip.1in

\proclaim{Claim 2.3} There is no infinite successive sequence
 of untwisting by the links under the case $\lambda > \mu$ with nonstationary
quasi-effective threshold.
\endproclaim

Suppose there is such an infinite sequence
$$
\CD
X = X_0@.\dashrightarrow @.X_1 @. \dasharrow @.X_2 @.\dashrightarrow @. @.\ \ \
\
\
\ \ \ \ \ \ \ \ \ \cdot\cdot\cdot
\ \ \ \ \ \ \ \ \ \ \ \ \ \dasharrow @.X_{k} @.\overset\psi_k\to\dashrightarrow
@.X_{k+1} @. \cdot\cdot\cdot \\
@V\phi VV @.@V\phi_1VV @. @V\phi_2VV @.@.\cdot\cdot\cdot @.
@V\phi_k VV @. @V\phi_{k+1}VV @. \cdot\cdot\cdot \\
S = S_0@.\leftarrow @. S_1 @.\leftarrow @.S_2 @.\leftarrow @.@.
\ \ \ \ \ \ \ \ \ \ \ \cdot\cdot\cdot
\ \ \ \ \ \ \ \ \ \ \ \ \ \ \ \ \leftarrow @.S_k @.\leftarrow @. S_{k+1} @.
\cdot\cdot\cdot
\\
\endCD
$$
Case: For some $k_0$, $\dim S_{k_0} \geq 1$.

\vskip.1in

In this case for $\forall k \geq k_0$ we have
$$\mu_k \in \frac{1}{3!}{\Bbb N}$$
as before, and $\{\mu_k\}$ is a nonstationary and nonincreasing infinte
sequence $\mu_0 \geq \mu_k > 0$, a contradiction!  In general, we only
need the boundedness of ${\Bbb Q}$-Fano $d$-folds for $d \leq n - 1$ up to this
point of the argument.

\vskip.1in

Finally

\vskip.1in

Case:For $\forall k$, $\dim S_k = 0$.

\vskip.1in

In this case, the $X_k$'s are all ${\Bbb Q}$-Fano variety with $\rho(X_k) =
1$.  Thus in dimension 3 Kawamata[8]'s result implies that they
belong to a bounded family.  Therefore, there exists $q \in {\Bbb N}$ s.t. $qD$
is Cartier for all integral Weil divisor on $X_k$ for $\forall k$ and there
exists $r \in {\Bbb N}$ s.t. $rK_{X_k}$ is Cartier for $\forall k$.  Then
another result of Kawamata[8] on the boundedness of the lengths of
the extremal rational curves says that there exists a rational curve $L_k$
on $X_k$ s.t. $0 < - K_{X_k} \cdot L_k \leq 2 \cdot \text{dim}X_k$, which
implies $$\mu_k \in \frac{1}{(r \cdot 2\text{dim}X)!q}{\Bbb N}.$$
Again $\{\mu_k\}$ is a nonstationary and  nonincreasing infinite sequence
$\mu_0 \geq \mu_k > 0$, a contradiction!

We remark that this last step is the only place where we use the boundedness
of ${\Bbb Q}$-Fano $n$-folds.

\vskip.1in

Cliams 2.1, 2.2 and 2.3 show that there is no infinite loop in the flowchart
of the Sarkisov program.

\vskip.1in

This completes the discussion of termination of the flowchart.

\newpage

\subhead{\S 3. Log Sarkisov Program with KLT Singularities}\endsubhead

\vskip.1in

In this section, we try to establish the Log Sarkisov Program for untwisting a
birational map between two log Mori fiber spaces with only kawamata log
terminal
singularities
$$
\CD
(X,B_X) @.\overset\Phi\to\dashrightarrow @.(X',B_{X'}) \\
@V\phi VV @.@VV\phi' V \\
S @.@. S' \\
\endCD
$$
The guiding principle throughout is that while the genuine Sarkisov program is
the one to untwist a birational map between two Mori fiber spaces obtained as
two different end results of $K$-MMP starting from one nonsingular projective
variety $W$
$$
\CD
@.@.W @.@.\\
@.K-\text{MMP}\swarrow@.@.\searrow K-\text{MMP}\\
X @.@.\overset\Phi\to\dashrightarrow @.@.X' \\
@V\phi VV @.@.@.@VV\phi' V \\
S @.@.@.@. S', \\
\endCD
$$
the log Sarkisov program should be the one to untwist a birational map between
two log Mori fiber spaces obtained as two different end results of $K +
B$-MMP starting from one log variety $(W,B_W)$ consisting of a nonsingular
projective variety $W$ and an S.N.C. divisor $B_W$
$$
\CD
@.@.(W,B_W) @.@.\\
@.K+B-\text{MMP}\swarrow@.@.\searrow K+B-\text{MMP}\\
(X,B_X) @.@.\overset\Phi\to\dashrightarrow @.@.(X',B_{X'}) \\
@V\phi VV @.@.@.@VV\phi' V \\
S @.@.@.@. S'. \\
\endCD
$$
(Note that in the two diagrams above we do not require a priori the slanted
arrows to be morphisms.)

While this principle does not put any restriction on the birational map $\Phi$
in the case of the genuine Sarkisov program (namely, for any birational map
$\Phi$ between two Mori fiber spaces we can find a common resolution $W$ s.t.
$X \rightarrow S$ and $X' \rightarrow S'$ are two end results of $K$-MMP as
the diagram above), this principle in the case of the log Sarkisov program
allows us to consider only such birational map $\Phi$ for which a log variety
as above exists.  This naturally leads to the notion of the Sarkisov relation:
(Two or more) Log Mori fiber spaces are Sarkisov related iff they are the end
results of
$K + B$-MMP starting from one log pair consisting of a nonsingular
projective variety and an S.N.C. divisor as a boundary.  Thus the log Sarkisov
program untwists a birational map between two log Mori fiber spaces which are
Sarkisov related, by factorizing it into links among intermediate Mori fiber
spaces all of which including the original two we require to be also Sarkisov
related.  (Remark again that in the case of the genuine Sarkisov program the
Sarkisov relation happens to be an equivalence relation and coincide with the
usual birational equivalence.)

Once we understand what the appropriate logarithmic generalization of
the Sarkisov program should be through the notion of the Sarkisov relation, the
flowchart for the log Sarkisov program with klt singularities works almost
parallel to that of the genuine Sarkisov program as well as termination except
the very last step.  In order to show that there is no infinite successive
sequence of links of type (II) with nonstationary quasi-effective threshold, we
have to use the conjecture by Borisov[3] in dimension 3
(cf.Nikulin[26]Alexeev[2]).

\proclaim{Conjecture 3.1} Fix a nonnegative rational number $0 \leq \epsilon
< 1$.  Then the family of log ${\Bbb Q}$-Fano $n$-folds (a normal projective
$n$-fold $X$ with only ${\Bbb Q}$-factorial log terminal singularities (and
thus automatically has only klt singularities) s.t. the anti-canonical
divisor $-K_X$ is ample) whose discrepancies are all $> - \epsilon$ is bounded.
\endproclaim

In dimension 2, the conjecture holds (cf.Nikulin[26]Alexeev[2]) and thus we
establish the log Sarkisov program for klt surfaces.  The conjecture in
dimension 3 for the case of Picard number 1 and $\epsilon = 0$ is the
theorem of Kawamata, which completes the proof of termination for the
genuine Sarkisov program for 3-folds.

\vskip.1in

In the following we discuss the log Sarkisov program with klt singularities
in detail.

\vskip.1in

\proclaim{Definition 3.2 (cf.Koll\'ar et al[16])} A log pair $(X,B_X = \Sigma
b_iB_i)$ has only kawamata log terminal singularities iff every discrete
valuation of $k(X)$ having center on $X$ has positive log discrepancy, i.e.,

o) $X$ is normal,

i) $0 \leq b_i < 1$ for $\forall i$, and

ii) there exists a log resolution $f:Y \rightarrow X$ where all the
$f$-exceptional divisors have positive log discrepancies (and thus this holds
for any log resolution).
\endproclaim

\vskip.1in

\proclaim{Definition 3.3} A log Mori fiber space $\phi:(X,B_X) \rightarrow
S$ with only klt singularities is the contraction of an extremal ray with
respect to $K_X + B_X$ from a log pair $(X,B_X)$ consisting of a normal
projective variety $X$ and a divisor $B_X$ with only ${\Bbb Q}$-factorial klt
singularities onto a variety $S$, i.e., $\phi$ is a morphism from a ${\Bbb
Q}$-factorial klt log pair $(X,B_X)$ with connected fibers onto a normal
variety
$S$ with
$\dim S < \dim X$ s.t. $\rho(X/S) = 1$ and $-(K_X + B_X)$ is
$\phi$-ample.
\endproclaim

\vskip.1in

\proclaim{Definition 3.4 (Sarkisov Relation)} Log Mori fiber spaces (resp. Log
minimal models) $$(X_0,B_{X_o}), (X_1,B_{X_1}), \cdot\cdot\cdot,(X_k,B_{X_k}),
\cdot\cdot\cdot,(X_l,B_{X_l})$$
are Sarkisov related iff they are all end results of $K + B$-MMP starting from
one appropriate log pair $(W,B_W)$ consisting of a nonsingular projective
variety
$W$ and an S.N.C. divisor $B_W = \Sigma b_iB_i$ with $0 \leq b_i \leq 1$.

A birational map $\Phi$ between two log Mori fiber spaces (resp. log minimal
models) which are Sarkisov related is by definition the one for which we have a
commutative diagram $$
\CD
@.@.(W,B_W) @.@.\\
@.p\swarrow@.@.\searrow q\\
(X,B_X) @.@.\overset\Phi\to\dashrightarrow @.@.(X',B_{X'}) \\
\endCD
$$
where $(W,B_W)$ is a log pair specified as above (Note that $p$ or $q$ may not
be a morphism.)
\endproclaim

\vskip.1in

The Sarkisov relation behaves very well for log Mori fiber spaces (or log
minimal models) with only kawamata log terminal singularities.

\vskip.1in

\proclaim{Proposition 3.5} Let
$$(X_0,B_{X_o}), (X_1,B_{X_1}), \cdot\cdot\cdot,(X_k,B_{X_k}),
\cdot\cdot\cdot,(X_l,B_{X_l})$$
be log Mori fiber spaces (resp. log minimal models) with only klt singularities
and $0 \leq \epsilon < 1$ a rational number such that all the coefficients of
the
boundaries $B_{X_k}$ are $\leq \epsilon$ and that all the discrepancies of the
log pairs $(X_k,B_{X_k})$ are $> - \epsilon$.

Then the following are equivalent:

\ \ (i) The log Mori fiber spaces (resp. log minimla models) with klt
singularities $$(X_0,B_{X_o}), (X_1,B_{X_1}), \cdot\cdot\cdot,(X_k,B_{X_k}),
\cdot\cdot\cdot,(X_l,B_{X_l})$$
are Sarkisov related, i.e., there exists a log variety $(W,B_W)$ consisting of
a nonsingular projective variety $W$ and an S.N.C. divisor $B_W$ as a boundary
such that all the log Mori fiber spaces (resp. log minimal models) are end
results of $K + B$-MMP over $Spec\ k$ starting from $(W,B_W)$.

\ \ (ii) There exists a log variety $(W,B_W)$ consisting of
a nonsingular projective variety $W$ and an S.N.C. divisor $B_W$ as a boundary
such that all the log Mori fiber spaces (resp. log minimal models) are end
results of $K + B$-MMP over $Spec\ k$ starting from $(W,B_W)$ and that
$$B_W = D_W(B_{X_0},B_{X_1}, \cdot\cdot\cdot,B_{X_k},
\cdot\cdot\cdot,B_{X_l}) + \Sigma_{E_j \text{not\ appearing\ as\ a\ divisor\
on\ any\ of\ }X_k}\epsilon E_j,$$
where
$$D_W(B_{X_0},B_{X_1}, \cdot\cdot\cdot,B_{X_k},
\cdot\cdot\cdot,B_{X_l}) := \Sigma d_mD_m$$
summation being taken over all divisors on $W$ which appear as a divisor on
some $X_m$ and $d_m$ being the coefficient of $D_m$ in $B_{X_m}$.

\ \ (iii) There exists a log variety $(W,B_W)$ consisting of
a nonsingular projective variety $W$ and an S.N.C. divisor $B_W$ as a boundary
such that each log Mori fiber space (resp. log minimal model) is dominated by a
birational morphism $p_k:(W,B_W) \rightarrow (X_k,B_{X_k})$ and an end result
of
$K + B$-MMP over $X_k$ starting from $(W,B_W)$ and that
$$B_W = D_W(B_{X_0},B_{X_1}, \cdot\cdot\cdot,B_{X_k},
\cdot\cdot\cdot,B_{X_l}) + \Sigma_{E_j \text{not\ appearing\ as\ a\ divisor\
on\ any\ of\ }X_k}\epsilon E_j.$$

\endproclaim

\demo{Proof}\enddemo The proposition is a straightforward consequence of the
following lemma, whose first claim holds not only for klt singularities but
also weakly kawamata log terminal singularities (or even more generally for log
canonical singularities) while the second claim only holds for klt
singularities.  This is why the Sarkisov relation behaves very well for klt
singularities but becomes quite subtle for wklt or lc singularities.  The
verification of the lemma is left to the reader as an exercise.

\proclaim{Lemma 3.6}

\ \ (i) Let $p:(W,B_W) \rightarrow (X,B_X)$ be a projective birational morphism
between ${\Bbb Q}$ factorial log varieties with klt (or more generally with
wklt or lc) singularities.  Then $p$ is a process of $K + B$-MMP over $X$
starting from $(W,B_W)$ iff $B_X = p_*(B_W)$ and the ramification divisor $R$
$$K_W + B_W = p^*(K_X + B_X) + R$$
has the same support as the exceptional locus $E(p)$ of $p$.

\ \ (ii) Let $(X,B_X)$ be a log pair with only ${\Bbb Q}$-factorial klt
singularities and $0 \leq \epsilon < 1$ a rational number such that all the
coefficients of the boundary $B_X$ are $\leq \epsilon$ and that all the
discrepancies are $> - \epsilon$.  Then any projective birational morphism
$p:(W,B_W) \rightarrow (X,B_X)$ from a log pair $(W,B_W)$ consisting of a
nonsingular projective variety and an S.N.C. divisor $B_W$
$$B_W = D_W(B_X) + \Sigma_{E_j \text{not\ appearing\ as\ a\ divisor\
on\ }X}\epsilon E_j,$$
is a process of $K + B$-MMP over $X$ starting from $(W,B_W)$.

\endproclaim

\proclaim{Remark 3.7}\endproclaim

Though unfortunately the Sarkisov relation is NOT an equivalence relation in
general, it is an equivalence relation for the following special classes of
log pairs $(X,B_X)$ with klt singularities: We fix $0 \leq \epsilon < 1$.
The class consists of ${\Bbb Q}$-factorial projective log pairs $(X,B_X)$
with klt singularities whose coefficients of the boundaries $B_X$ are all equal
to $\epsilon$, and all the discrepancies of the valuations of the exceptional
divisors are $> - \epsilon$.

The genuine Sarkisov program is nothing but the program for the class given
by $\epsilon = 0$.

\vskip.1in

\proclaim{Theorem 3.8 (Log Sarkisov Program for log 3-folds with klt
singularities)} Let  $$
\CD
(X,B_X) @.\overset\Phi\to\dashrightarrow @.(X',B_{X'}) \\
@V\phi VV @.@VV\phi' V \\
S @.@. S' \\
\endCD
$$
be a birational map between two log Mori fiber spaces in dimension 3 with only
klt singularities, which are Sarkisov related.

Suppose the Borisov conjecture holds for ${\Bbb Q}$-Fano 3-folds with klt
singularities and Picard number 1.

Then for any rational number $0 \leq \epsilon < 1$ such
that all the coefficients in $B_X$ or $B_{X'}$ are $\leq \epsilon$ and that all
the discrepancies of $(X,B_X)$ or $(X',B_{X'})$ are $> - \epsilon$, $\Phi$ is a
composite of 4 types of links as in the genuine Sarkisov program

$$
\CD
(X,B_X) = (X_0,B_{X_0}) @.\dashrightarrow @. (X_1,B_{X_1}) @. \dasharrow @.
\cdot\cdot\cdot @. \dashrightarrow @. (X_k,B_{X_k}) @. \dashrightarrow
@. \cdot\cdot\cdot @.\dashrightarrow @. (X',B_{X'})\\
@V\phi VV @. @V\phi_1VV
@. @. \cdot\cdot\cdot @. @V\phi_kVV @.@.\cdot\cdot\cdot @. @V\phi'VV \\
S = S_0 @.@. S_1 @.@.\cdot\cdot\cdot @.@. S_k @.@.\cdot\cdot\cdot @.@.S' \\
\endCD
$$

such that all the log Mori fiber spaces $(X_k,B_{X_k})$ have only ${\Bbb
Q}$-factorial klt singularities, the coefficients of $B_{X_k}$ are $\leq
\epsilon$, all the discrepancies of $(X_k,B_{X_k})$ are $> - \epsilon$ and
that all the $(X_k,B_{X_k})$ are Sarkisov related.  More precisely, all the
log Mori fiber spaces are dominated by birational morphisms
$$p_k:(W,B_W) \rightarrow (X_k,B_{X_k})$$
from a log pair $(W,B_W)$ consisting of a nonsingular projective 3-fold and
an S.N.C. divisor $B_W$
$$B_W = D_W(B_X,B_{X'}) + \Sigma_{E_j \text{not\ appearing\ as\ a\ divisor\
either\ on\ }X \text{\ or\ on\ }X'}\epsilon E_j
$$
and each $(X_k,B_{X_k})$ is an end result of $K + B$-MMP over $X_k$ starting
from $(W,B_W)$.
\endproclaim

The strategy to establish the log Sarkisov program goes along the same line as
the one to establish the genuine Sarkisov program, constructing $(W,B_W)$ as
above inductively as the program proceeds.

\vskip.1in

We define the log Sarkisov degree of an intermediate log Mori fiber space
$(X_k,B_{X_k})$ which appears in the due course of untwisting a birational map
$$
\CD
(X,B_X) @.\overset\Phi\to\dashrightarrow @.(X',B_{X'}) \\
@V\phi VV @.@VV\phi' V \\
S @.@. S' \\
\endCD
$$
between two log Mori fiber spaces which are Sarkisov related as follows.

\vskip.1in

\proclaim{Definition 3.9 (the log Sarkisov degree)} Let
$$
\CD
(X,B_X) @.\overset\Phi\to\dashrightarrow @.(X',B_{X'}) \\
@V\phi VV @.@VV\phi' V \\
S @.@. S' \\
\endCD
$$
be a birational map between two log Mori fiber spaces with only klt
singularities which are Sarkisov related, and fix a rational
number $0 \leq \epsilon < 1$ such that all the
coefficients in $B_X$ or $B_{X'}$ are $\leq \epsilon$ and that all the
discrepancies of $(X,B_X)$ or $(X',B_{X'})$ are $> - \epsilon$.  Then the log
Sarkisov degree of any intermediate log Mori fiber space
$\phi_k:(X_k,B_{X_k}) \rightarrow S_k$ that appears in the due course of
untwisting the birational map with reference to the fixed log Mori fiber space
$\phi':(X',B_{X'})
\rightarrow S'$ is the triplet $$(\mu_k,\lambda_{\epsilon k},e_{\epsilon k})$$
of
the numbers defined below, endowed with the lexicographical order.
\endproclaim

Notice that there is an auxiliary parameter
$\epsilon$, which was implicit (actually equal to 0) in the case of the genuine
Sarkisov degree.  Also note that the log Sarkisov degree depends not only on
$(X_k,B_{X_k})$ and $(X',B_{X'})$ but also on the initial log Mori fiber
space $(X,B_X)$.

\vskip.1in

First we take a very ample divisor $A'$ on $S'$ and a sufficiently
divisible $\mu' \in {\Bbb N}$ such that
$${\Cal H}_{X'} = -\mu'(K_{X'} + B_{X'}) + {\phi'}^*A'$$
is very ample on $X'$.  ${\Cal H}_{X_k}$ is the strict transform of ${\Cal
H}_{X'}$ on $X_k$.

\vskip.1in

$\boxed{\mu_k:\text{the\ quasi-effective\ threshold}}$

\vskip.1in

The quasi-effective threshold is defined exactly the same way as before,
replacing $K$ with $K + B$, namely
$$\mu_k \in {\Bbb Q}_{>0}\ \mu_k(K_{X_k} + B_{X_k}) + {\Cal H}_{X_k} \equiv 0
\text{\ over\ } S_k$$

\vskip.1in

$\boxed{\lambda_{\epsilon k}:\text{the\ maximal\ multiplicity\ of\ an\
extremal\
ray}}$

\vskip.1in

Let $(W,B_W)$ be a log pair consisting of a nonsingular projective variety $W$
and an S.N.C. divisor $B_W$
$$B_W = D_W(B_X,B_{X'}) + \Sigma_{E_j \text{not\ appearing\ as\ a\ divisor\
either\ on\ }X \text{\ or\ on\ }X'}\epsilon E_j
$$
such that it dominates
$$\align
p:(W,B_W) &\rightarrow (X,B_X)\\
q:(W,B_W) &\rightarrow (X',B_{X'})\\
p_k:(W,B_W) &\rightarrow (X_k,B_{X_k})\\
\endalign$$
by birational morphisms and that they are all processes of $K + B$-MMP over
$X, X'$ and $X_k$, respectively.  (The existence of such $(W,B_W)$ will be
shown inductively in the course of log Sarkisov program.)

Then
$$\frac{1}{\lambda_{\epsilon k}} := \text{max}\{c \in {\Bbb Q}_{> 0};(K_W +
B_W)
+ c{\Cal H}_{X_k} = p_k^*(K_{X_k} + B_{X_k}) + c{\Cal H}_{X_k} + \text{some\
effective\ divisor}\}.$$
Note that $\frac{1}{\lambda_{\epsilon k}}$ is independent of the choice of
such $(W,B_W)$ and well-defined.  When $Bs({\Cal H}_X) = \emptyset$,
$\lambda_{\epsilon k} = 0$ by definition.

We note that in general
$$\lambda_{\epsilon k} \leq \epsilon-\text{log\ terminal\ threshold\ of\
}(X_k,B_{X_k}) \text{\ w.r.t.\ }{\Cal H}_{X_k}.$$
The $\epsilon$-log termianl threshold of the pair $(X_k,B_{X_k})$ is defined
to be
$$\text{max}\{c \in {\Bbb Q}_{> 0};(K_{V_k} + B_{V_k})
+ c{\Cal H}_{V_k} = v_k^*((K_{X_k} + B_{X_k}) + c{\Cal H}_{X_k}) + \text{some\
effective\ divisor}\},$$
where $V_k$ is any nosingular projective variety which dominates both
$v_k:V_k \rightarrow X_k$ and $X'$
 by birational morphisms such that the union of the exceptional locus of $v_k$
and ${v_k}_*^{-1}(B_{X_k})$ is an S.N.C. divisor (Recall the ${\Bbb
Q}$-factoriality of $X_k$), and
$$B_{V_k} = {v_k}^{-1}_*(B_{X_k}) + \Sigma_{E_j \text{not\ Apearing\ as\ a\
divisor\ on\ }X_k}\epsilon E_j.$$

\vskip.1in

$\boxed{e_{\epsilon k}:\text{the\ number\ of\ }K + B
+ \frac{1}{\lambda_{\epsilon k}}{\Cal H}-\text{crepant\ divisors}}$

\vskip.1in

$e_{\epsilon k}$ is defined to be the number of exceptional divisors for $p_k$
whose coefficient in the ramification divisor $R$ is $0$
$$(K_W + B_W) + \frac{1}{\lambda_{\epsilon k}}{\Cal H}_W = p_k^*((K_{X_k} +
B_{X_k}) + \frac{1}{\lambda_{\epsilon k}}{\Cal H}_{X_k}) + R.$$
Again this is independent of the choice of $(W,B_W)$ and well-defined.

\vskip.1in

\proclaim{N\"other-Fano Criterion for the log Sarkisov program with KLT
singularities} The birational map $\Phi_k$ between an intermediate log Mori
fiber space $\phi_k:(X_k,B_{X_k}) \rightarrow S_k$ and $\phi':(X',B_{X'})
\rightarrow S'$ is an isomorphism of log Mori fiber spaces
$$
\CD
(X_k,B_{X_k}) @.\overset\Phi_k\to{\overset\sim\to\rightarrow} @.(X',B_{X'}) \\
@V\phi_k VV @.@VV\phi' V \\
S_k @.\overset\sim\to\rightarrow @. S' \\
\endCD
$$
if $\lambda_{\epsilon k} \leq \mu_k$ and $(K_{X_k} + B_{X_k}) +
\frac{1}{\lambda_{\epsilon k}}{\Cal H}_{X_k}$ is nef.
\endproclaim

The proof goes verbatim to the one for the genuine Sarkisov program, taking
the common resolution $(W,B_W)$ which Sarkisov-relates $(X_k,B_{X_k})$ and
$(X',B_{X'})$ into consideration.

\vskip.2in

$$\bold{Flowchart\ for\ Log\ Sarkisov\ Program\ with\ KLT\ Singularities}$$

The flowchart for the log Sarkisov program with klt singularities goes almost
parallel to the one for the genuine Sarkisov program replacing $K$ with
$K + B$, constructing inductively a log pair $(W,B_W)$ which dominates all the
intermediate log Mori fiber spaces in the process of untwisting.

\vskip.1in

Let
$$
\CD (X,B_X) @.\overset\Phi\to\dashrightarrow @.(X',B_{X'}) \\ @V\phi VV
@.@VV\phi' V \\ S @.@. S' \\
\endCD
$$
be a birational map between two log Mori fiber spaces with only klt
singularities, which are Sarkisov related.  We fix a rational number $0 \leq
\epsilon < 1$ such that all the coefficients in $B_X$ or $B_{X'}$ are
$\leq \epsilon$ and that all the discrepancies of $(X,B_X)$ or $(X',B_{X'})$
are
$> - \epsilon$.

We remark that for all the relevant log pairs $(U,B_U)$ that appear in the
course of the log Sarkisov program (including the auxiliary resolutions we
take) the boundaries $B_U$ are always taken to be of the form
$$B_U = D_U(B_X,B_{X'}) + \Sigma_{E_j \text{\ not\ appearing\ as\ a\ divisor\
either\ on\ }X \text{\ or\ on\ }X'}\epsilon E_j.$$

Before we start the flowchart, we note that by Proposition 3.5 there exists a
log variety $(W_0,B_{W_0})$ consisting of a nonsingular projective variety
$W_0$ and an S.N.C. divisor $B_{W_0}$ as a boundary such that $(X,B_X) =
(X_0,B_{X_0})$ (resp.
$(X',B_{X'})$) is dominated by a birational morphism
$p = p_{0,0}:(W_0,B_{W_0}) \rightarrow (X_0,B_{X_0})$ (resp. $q =
q_0:(W_0,B_{W_0})
\rightarrow (X',B_{X'})$) and an end result of
$K + B$-MMP over
$X_0$ (resp. over $X'$) starting from $(W_0,B_{W_0})$ and that
$$B_{W_0} = D_{W_0}(B_X,B_{X'}) + \Sigma_{E_j \text{not\ appearing\ as\ a\
divisor\ on\ either\ }X \text{\ or\ on\ }X'}\epsilon E_j.$$

\vskip.1in

Suppose we have untwisted the birational map up to the $k$-th stage and
constructed a log pair $(W_k,B_{W_k})$ consisting of a nonsingular projective
variety $W_k$ and an S.N.C. divisor $B_{W_k}$ as a boundary such that
each log Mori fiber space $(X_m,B_{X_m})\ m = 0, 1, \cdot\cdot\cdot, k$ (and
$(X',B_{X'})$) is dominated by a birational morphism
$p_{k,m}:(W_k,B_{W_k}) \rightarrow (X_m,B_{X_m})$ (and
$q_k:(W_k,B_{W_k})
\rightarrow (X',B_{X'})$) and an end result of
$K + B$-MMP over
$X_m$ (and over $X'$) starting from $(W_k,B_{W_k})$ and that
$$B_{W_k} = D_{W_k}(B_X,B_{X'}) + \Sigma_{E_j \text{not\ appearing\ as\
a\ divisor\ on\ either\ }X \text{\ or\ on\ }X'}\epsilon E_j.$$

\vskip.1in

The first question then to ask as in the genuine Sarkisov program is:
$$\lambda_{\epsilon k} > \mu_k?$$

$$\boxed{\text{Case:}\lambda_{\epsilon k} \leq \mu_k}$$

\vskip.1in

In this case, thanks to the N\"other-Fano inequality for klt singularities,
the program works completely parallel to the genuine Sarkisov program.  If
$(K_{X_k} + B_{X_k}) + \frac{1}{\lambda_{\epsilon k}}{\Cal H}_{X_k}$ is
nef, then the program comes to an end.  If not, then after untwisting the
birational map either by a link of type (III) or (IV), the quasi-effective
threshold strictly drops $$\mu_{k+1} < \mu_k$$ as before.

We also take a nonsingular projective variety $W_{k+1}$ by blowing up $W_k$
further so that
it dominates each log Mori fiber space by a birational morphism
$p_{k+1,m}:(W_{k+1},B_{W_{k+1}}) \rightarrow (X_m,B_{X_m})\ m = 0,1,
\cdot\cdot\cdot, k, k+1$ (and
$q_{k+1}:(W_{k+1},B_{W_{k+1}}) \rightarrow (X',B_{X'})$) and that
$$B_{W_{k+1}} := D_{W_{k+1}}(B_X,B_{X'}) + \Sigma_{E_j \text{not\ appearing\
as\ a\ divisor\
on\ either\ }X \text{\ or\ on\ }X'}\epsilon E_j$$
is an S.N.C. divisor.  Then by Lemma 3.6 $(W_k,B_{W_k})$ is an end result of $K
+
B$-MMP over
$W_k$ starting from $(W_{k+1},B_{W_{k+1}})$ and thus each log Mori fiber space
$(X_m,B_{X_m}) m = 0, 1, \cdot\cdot\cdot, k$ (and $(X',B_{X'})$) is an end
result of $K +
B$-MMP starting from $(W_{k+1},B_{W_{k+1}})$ over $X_m$.  As for
$(X_{k+1},B_{X_{k+1}})$
which is an end result of $K + B + \frac{1}{\lambda_{\epsilon k}}{\Cal H}$-MMP
starting
from $(X_k,B_{X_k})$ (over
$T$), all the coeffiecients of $B_{X_{k+1}}$ are $\leq \epsilon$ by
construction.
Moreover, for any valuation $E$ whose center on $X_{k+1}$ has codimension at
least 2
$$a(E,X_{k+1},B_{X_{k+1}}) = a(E,X_k,B_{X_k}) < - \epsilon$$
if the birational map $(X_k,B_{X_k}) \dashrightarrow (X_{k+1},B_{X_{k+1}})$ is
isomorphic
at the center of the valuation $E$, and
$$\align
a(E,X_{k+1},B_{X_{k+1}}) &\leq a(E,X_{k+1},B_{X_{k+1}} +
\frac{1}{\lambda_{\epsilon k}}{\Cal
H}_{X_{k+1}})\\
 &< a(E,X_k,B_{X_k} + \frac{1}{\lambda_{\epsilon k}}{\Cal H}_{X_{k}}) \leq -
\epsilon.\\
\endalign$$
if the birational map $(X_k,B_{X_k}) \dashrightarrow (X_{k+1},B_{X_{k+1}})$ is
not isomorphic
at the center of the valuation $E$.

Thus again by Lemma 3.6 $(X_{k+1},B_{X_{k+1}})$ is an end result of $K + B$-MMP
starting from $(W_{k+1},B_{W_{k+1}})$ over $X_{k+1}$.  Therefore,
$(W_{k+1},B_{W_{k+1}})$ satisfies the desired inductive property.

\vskip.1in

$$\boxed{\text{Case:}\lambda_{\epsilon k} > \mu_k}$$

\vskip.1in

In this case we
$$\text{Take\ a\ maximal\ divisorial\ blow\ up\ }d_k:(Z_k,B_{Z_k}) \rightarrow
(X_k,B_{X_k})$$
with respect to $(K_{X_k} + B_{X_k}) + \frac{1}{\lambda_{\epsilon k}}{\Cal
H}_{X_k}$, and in
the log Sarkisov program we also require $(Z_k,B_{Z_k})$ is obtained through
$K + B + \frac{1}{\lambda_{\epsilon k}}{\Cal H}$-MMP possibly followed by $K +
B$-MMP starting from
$(W_k,B_{W_k})$ over $X_k$, i.e.,

\ \ o) in the proof of Proposition 1.5 instead of carrying out a $K +
\frac{1}{\lambda}{\Cal H}$-MMP over $X$ starting from $Y$ to get a minimal
model
$(Z',\frac{1}{\lambda}{\Cal H})$ and then running a $K$-MMP over $X$ to obtain
a maximal
divisorial blow up, we carry out a
$K + B + \frac{1}{\lambda_{\epsilon k}}{\Cal H}$-MMP over $X_k$ starting from
$(W_k,B_{W_k})$ to get a minimal model $({Z'}_k + B_{{Z'}_k} +
\frac{1}{\lambda_{\epsilon
k}}{\Cal H}_{{Z'}_k})$ and then run a $K + B$-MMP over $X_k$ to obtain the
maximal divisorial
blow up $d_k:(Z_k,B_{Z_k}) \rightarrow (X_k,B_{X_k})$,

\ \ i) $\rho(Z_k/X_k) = 1$,

\ \ ii) the exceptional locus of $d_k$ is a prime divisor $E_k$, and

\ \ iii) $d_k$ is $K + B + \frac{1}{\lambda_{\epsilon k}}{\Cal H}$-crepant,
i.e.,
$$(K_Z + B_Z) + \frac{1}{\lambda_{\epsilon k}}{\Cal H}_Z = d_k^*((K_{X_k} +
B_{X_k}) +
\frac{1}{\lambda_{\epsilon k}}{\Cal H}_{X_k}).$$

Then the rest goes parallel to the genuine Sarkisov program, untwisting the
birational map
further by a link of either type (II) or (I).  After untwisting, we have
$$\mu_{k+1} \leq \mu_k$$
and if $\mu_{k+1} = \mu_k$ and $\dim S_k = \dim S_{k+1}$, then $\psi_k$ is a
square.  Moreover,
$$\lambda_{\epsilon, k+1} \leq \lambda_{\epsilon k}$$
and if $\lambda_{\epsilon, k+1} = \lambda_{\epsilon k}$ then $e_{\epsilon, k+1}
<
e_{\epsilon k}$.

We also take a nonsingular projective variety $W_{k+1}$ by blowing up $W_k$
further so that it
dominates each log Mori fiber space by a birational morphism
$p_{k+1,m}:(W_{k+1},B_{W_{k+1}}) \rightarrow (X_m,B_{X_m})\ m = 0, 1,
\cdot\cdot\cdot, k, k+1$ (and
$q_{k+1}:(W_{k+1},B_{W_{k+1}}) \rightarrow (X',B_{X'})$) and that
$$B_{W_{k+1}} := D_{W_{k+1}}(B_X,B_{X'}) + \Sigma_{E_j \text{not\ appearing\
as\ a\ divisor\
on\ either\ }X \text{\ or\ on\ }X'}\epsilon E_j$$  is an S.N.C. divisor.  Then
by Lemma
3.6
$(W_k,B_{W_k})$ is an end result of $K + B$-MMP over
$W_k$ starting from $(W_{k+1},B_{W_{k+1}})$ and thus each log Mori fiber space
$(X_m,B_{X_m})\ m = 0, 1, \cdot\cdot\cdot, k$ (and $(X',B_{X'})$) is an end
result of $K + B$-MMP starting from $(W_{k+1},B_{W_{k+1}})$ over $X_m$.

As for $(X_{k+1},B_{X_{k+1}})$, note that it is an end result of $K + B +
\frac{1}{\lambda_{\epsilon k}}{\Cal H}$-MMP starting from
$(Z_k,B_{Z_k})$ (over
$S_k$), which itself is an end result of $K + B + \frac{1}{\lambda_{\epsilon
k}}{\Cal
H}$-MMP possibly followed by $K + B$-MMP over $X_k$ starting from
$(W_k,B_{W_k})$.  All the
coeffiecients of $B_{X_{k+1}}$ are $\leq \epsilon$ by construction.

Moreover, for any valuation $E$ whose center on ${Z'}_k$ has at least
codimension 2
$$a(E,{Z'}_k,B_{{Z'}_k}) = a(E,W_k,B_{W_k}) < - \epsilon$$
if the birational map
$(W_k,B_{W_k}) \dashrightarrow ({Z'}_k,B_{{Z'}_k})$ is isomorphic at the center
of the
valuation $E$, and
$$\align
a(E,{Z'}_k,B_{{Z'}_k}) &\leq a(E,{Z'}_k,B_{{Z'}_k} + \frac{1}{\lambda_{\epsilon
k}}{\Cal H}_{{Z'}_k})\\
 &< a(E,W_k,B_{W_k} + \frac{1}{\lambda_{\epsilon k}}{\Cal H}_{W_{k}}) \leq -
\epsilon\\
\endalign$$
if the birational map $(W_k,B_{W_k}) \dashrightarrow ({Z'}_k,B_{{Z'}_k})$ is
not isomorphic at the center of the valuation $E$.

For any valuation $E$ whose center on $Z_k$ has at least codimension 2
$$a(E,Z_k,B_{Z_k}) = a(E,{Z'}_k,B_{{Z'}_k}) < - \epsilon$$
if the birational map
$({Z'}_k,B_{{Z'}_k}) \dashrightarrow (Z_k,B_{Z_k})$ is isomorphic at the center
of the
valuation $E$, and
$$a(E,Z_k,B_{Z_k}) < a(E,{Z'}_k,B_{{Z'}_k}) \leq - \epsilon$$
if the birational map $({Z'}_k,B_{{Z'}_k}) \dashrightarrow (Z_k,B_{Z_k})$ is
not
isomorphic at the center of the valuation $E$.

Finally for any valuation $E$ whose center on $X_{k+1}$ has at least
codimension 2
$$a(E,X_{k+1},B_{X_{k+1}}) = a(E,Z_k,B_{Z_k}) < - \epsilon$$
if the birational map
$(Z_k,B_{Z_k}) \dashrightarrow (X_{k+1},B_{X_{k+1}})$ is isomorphic at the
center of the
valuation $E$, and
$$\align
a(E,X_{k+1},B_{X_{k+1}}) &\leq a(E,X_{k+1},B_{X_{k+1}} +
\frac{1}{\lambda_{\epsilon k}}{\Cal H}_{X_{k+1}})\\
 &< a(E,Z_k,B_{Z_k} + \frac{1}{\lambda_{\epsilon k}}{\Cal H}_{Z_{k}}) \leq -
\epsilon\\
\endalign$$
if the birational map $(Z_k,B_{Z_k}) \dashrightarrow (X_{k+1},B_{X_{k+1}})$ is
not
isomorphic at the center of the valuation $E$.

Thus again by Lemma 3.6 $(X_{k+1},B_{X_{k+1}})$ is an end result of $K + B$-MMP
starting from
$(W_{k+1},B_{W_{k+1}})$ over $X_{k+1}$.  Therefore, $(W_{k+1},B_{W_{k+1}})$
satisfies the
desired inductive property.

\vskip.2in

$$\bold{Termination\ of\ Flowchart\ for\ Log\ Sarkisov\ Program\ with\ KLT\
Singularities}$$

\vskip.1in

One of the key points to show the terminatin of of the flowchart for the
genuine Sarkisov
program is the discreteness of the quasi-effective thresholds, which follows
from the
boundedness of ${\Bbb Q}$-Fano varieties which are fibers of the Mori fiber
spaces (a
nonsingular rational curve, Del Pezzo surfaces and ${\Bbb Q}$-Fano 3-folds
for the genuine Sarkisov program in dimension 3).  In the log Sarkisov program,
we rely on the Borisov conjecture, which is a theorem in dimension 2 thanks to
Nikulin[26][Alexeev[2], to show the discreteness of the quasi-effective
thresholds.

\vskip.1in

\proclaim{Claim 3.10} Let $\phi_k:(X_k,B_{X_k}) \rightarrow S_k$ be a log Mori
fiber space in dimension 3 with only klt singularities in the process of the
log
Sarkisov program untwisting a birational map
$$
\CD
(X,B_X) @.\overset\Phi\to\dashrightarrow @.(X',B_{X'}) \\
@V\phi VV @.@VV\phi' V \\
S @.@. S' \\
\endCD
$$
between two log Mori fiber spaces in dimension 3 with only klt singularities
which are Sarkisov related, fixing $\epsilon$ as before.

If $\dim S_k \geq 1$, then the denominator of the
quasi-effective threshold $\mu_k$ is universally bounded by a fixed constant
depending only on $\epsilon$ and the coefficients of $B_X$ and $B_{X'}$.
\endproclaim

\demo{proof}\enddemo

Let $d$ be the l.c.m. of the denominator of $\epsilon$ and those of
coefficients of $B_X$ and $B_{X'}$.

When $\dim S_k = 2$, $l$ being a nonsingular rational curve which is a general
fiber
of $\phi_k$,
$$\align
-2 \leq (K_{X_k} + B_{X_k}) \cdot l &< 0\\
(K_{X_k} + B_{X_k}) \cdot l &\in \frac{1}{d}{\Bbb Z}_{<0}\\
\{\mu_k(K_{X_k} + B_{X_k}) + {\Cal H}_{X_k}\} \cdot l  &= 0\\
\endalign$$
imply
$$\mu_k \in \frac{1}{(2d)!}{\Bbb N}.$$
When $\dim S_k = 1$, a general fiber $F_k$ is a log Del Pezzo surface (a normal
projective surface with only quotient singularities having an ample
anti-canonical divisor)
whose discrepancies are all $> - \epsilon$.  Therefore, by
Nikulin[26][Alexeev[2] we conclude that the family of such log Del Pezzo
surfaces
is bounded.  Now $q$ being the universal ${\Bbb Q}$-factorial index for such
surfaces and $r$ the index for the canonical divisors, we have
$$\mu_k \in \frac{1}{(r \cdot 2 \cdot \dim X)!q}{\Bbb N}.$$
In arbitrary dimension, we have to use the Borisov conjecture for boundedness
of log ${\Bbb Q}$-Fano $d$-folds for $d \leq n - 1$ in general to derive this
claim.

\vskip.1in

Thanks to Claim 3.10, the argument for termination goes parallel replacing
$K$ with $K + B$ for Claims 2.1, 2.2 and the first case in Claim 2.3 to that of
termination of the genuine Sarkisov program.  (We note that in Step 3 of the
proof of Claim 2.2 we repalce the local canonical threshold with the local
version of $\frac{1}{\lambda_{\epsilon}}$.)  We have to use a conjecture of
Borisov to establish the last step of the second case in Claim 2.3.

\proclaim{Conjecture 3.11 (cf.Borisov[3])} Fix a rational number
$0 \leq \epsilon < 1$.  Then the family of log ${\Bbb Q}$-Fano 3-folds
(normal projective 3-folds $X$ with only ${\Bbb Q}$-factorial log terminal
(equivalently, klt) singularities s.t. the anti-canonical divisors $-K_X$ are
ample) with Picard number 1 and whose discrepancies are all $> - \epsilon$,
is bounded.
\endproclaim

In general, to establish Claim 2.3 in the second case we need the conjecture
above for log ${\Bbb Q}$-Fano $n$-folds.

\vskip.1in

This completes the discussion of termination of log Sarkisov program with
klt singularities.

\newpage

\subhead{\S 4. Log Sarkisov Program with WKLT Singularities}\endsubhead

\vskip.1in

In this section, we establsih the log Sarkisov program with weakly kawamata log
terminal singularities in dimension 2, and then discuss the problems one has to
face attempting to establish the log Sarkisov program with weakly kawamata log
terminal singularities in higher dimension.  We also prove that Sarkisov
related log minimal 3-folds with klt or wklt singularities are connected by a
sequence of log flops.

\vskip.1in

\proclaim{Definition 4.1 (Local) (cf.Koll\'ar et al[16])} Let $(X,B_X)$ be a
germ (with respect to Zariski topology) around a point $P \in X$.  $(X,B_X)$
has only weakly kawamata log terminal singularities if there exists a Zariski
open set $P \in U \subset X$ s.t. there exists a log resolution $f:V
\rightarrow U$ such that all the log discrepancies of the exceptional
divisors with center on $U$ are positive and we have an $f$-anti-ample
effective divisor whose support coincides with that of the exceptional locus
of $f$.  (Note that $B_X$ may have components with coefficient 1.)
\endproclaim

The relation between local and global properties of wklt singularities was
clarified by the following result of Szab\'o[32].

\proclaim{Proposition 4.2 (Global) (cf.Szab\'o[32])} Let $(X,B_X)$ be a
projective log variety which has locally only weakly kawamata log terminal
singularities.  Then there exists a log resolution (global) $f:Y \rightarrow
X$ such that all the log discrepancies of the exceptional divisors with
center on $X$ are positive and we have an $f$-anti-ample effective divisor
whose support coincides with that of the exceptional locus of $f$.
\endproclaim

\proclaim{Corollary 4.3 (Characterization of a log Mori fiber space (or a log
minimal model) with wklt singularities)} A log Mori fiber space $\phi:(X,B_X)
\rightarrow S$ (resp. a log minimal model $(X,B_X)$) in dimension $n$ ($\leq
3$)
with only ${\Bbb Q}$-factorial wklt singularities is an end result of a $K +
B$-MMP starting from $(Y,B_Y)$ where $Y$ is a nonsingular projective n-fold and
$B_Y = \Sigma b_iB_i$ is an S.N.C. divisor with $0 \leq b_i \leq 1$, and the
converse holds, i.e., any end result of fibering type (resp. of minimal model
type) of a $K + B$-MMP starting from $(Y,B_Y)$ as above is a log Mori fiber
space
(resp. a log minimal model) with only ${\Bbb Q}$-factorial wklt singularities.
(Once we have the log-MMP in dimension $n$, the same statement holds in
dimension
$n$.)  \endproclaim

\vskip.1in

We prove the well-behavior of the Sarkisov relation for log minimal
models with wklt singularities as follows, thanks to the fact that the nef log
canonical divisors of the Sarkisov related log minimal models are all
essentially
the same and uniquely characterized as the nef part of the Zariski
decomposition
of the log canonical divisor of an arbitrary log resolution.  However, we fail
to
prove the well-behavior of the Sarkisov relation for log Mori fiber spaces
with wklt singularities in dimension $> 2$, mainly because the lack of the
Zariski decomposition and the failure for the statement (ii) of
Proposition 3.5 to hold.

\proclaim{Proposition 4.4 = Proposition 3.5 for log minimal models with wklt
singularities} Let
$$(X_0,B_{X_o}), (X_1,B_{X_1}), \cdot\cdot\cdot,(X_k,B_{X_k}),
\cdot\cdot\cdot,(X_l,B_{X_l})$$
be log minimal models with only ${\Bbb Q}$-factorial wklt singularities.  Then
(i) (ii) and (iii) as in Proposition 3.5 (replacing the assumption of klt
singularities with that of wklt singularities and allowing the possibility
$\epsilon = 1$) are equivalent.
\endproclaim

\demo{Proof}\enddemo The implications $(iii) \Rightarrow (ii) \Rightarrow
(i)$ are obvious (rgardless whether they are log minimal models or log Mori
fiber spaces).  We only have to prove $(i) \Rightarrow (iii)$.

Take a log pair $(W,B_W)$ as in (i) and let
$$p_k:(W,B_W) \dashrightarrow (X_k,B_{X_k})$$
be a birational map which is a $K + B$-MMP over $Spec k$.  We denote by $C_k$
the closed set in $W$ so that
$$p_k:(W - C_k,B_W|_{W - C_k}) \overset\sim\to\rightarrow (X_k -
I_k,B_{X_k}|_{X_k - I_k}),$$
where $I_k$ is the indeterminacy of the birational map ${p_k}^{-1}$.  We can
take a blowup $\sigma:W' \rightarrow W$ whose centers are all over $\cup_k C_k$
such that each $X_k$ is dominated by a birational morphism
${p_k}':W' \rightarrow X_k$ and that $\sigma^{-1}_*(B_W) \cup E(\sigma)$ is an
S.N.C. divisor.

Set
$${\hat B}_{W'} := \sigma^{-1}_*(B_W) + \Sigma_{E_j \text{not\ appearing\ as\
a\
divisor\ on\ }W}\epsilon E_j.$$

For each $k$ we have the ramification formulae
$$K_{W'} + {\hat B}_{W'} = \sigma^*(K_W + B_W) + R_{\sigma}$$
where $R_{\sigma}$ is an effective ramification divisor (whose support may not
coincide with $E(\sigma)$), and
$$\sigma^*(K_W + B_W) = \sigma^*\{p_k^*(K_{X_k} + B_{X_k}) +
R_{p_k}\}$$
where a $p_k \circ \sigma$-exceptional divisor $E$ has a
strictly positive coefficient in $\sigma^*R_{p_k}$ iff the center of $E$ on $W$
is contained in $C_k$.

Now an easy application of the Negativity Lemma (cf.Koll\'ar[15],Lemma 4.3)
shows that
$$(p_0 \circ \sigma)^*(K_{X_0} + B_{X_0}) = \cdot\cdot\cdot = (p_k \circ
\sigma)^*(K_{X_k} + B_{X_k}) =
\cdot\cdot\cdot = (p_l \circ \sigma)^*(K_{X_l} + B_{X_l})$$
giving the nef part of the Zariski decomposition of $K_{W'} + {\hat B}_{W'}$,
and hence
$$R_{\sigma} + \sigma^*R_{p_0} = \cdot\cdot\cdot = R_{\sigma} + \sigma^*R_{p_k}
= \cdot\cdot\cdot = R_{\sigma} +
\sigma^*R_{p_l}.$$
Now any $k$, a $p_k \circ \sigma$-exceptional divisor $E$ has
the center on $W$ contained in $\cup_i I_i$ and thus has a strictly positive
coefficient in $\sigma^*R_{p_i}$ for some $i$.  Therefore, the above equality
implies it has a strictly positive coefficient in $R_{\sigma} +
\sigma^*R_{p_k}$.  Thus
$$K_{W'} + {\hat B}_{W'} = (p_k \circ \sigma)^*(K_{X_k} + B_{X_k}) + {\hat
R}_{p_k \circ \sigma}$$
where
$$\text{supp}\ {\hat
R}_{p_k \circ \sigma} = \text{supp}\ E(p_k \circ \sigma).$$
Therefore, finally by setting
$$B_{W'} = D_W(B_{X_0},B_{X_1}, \cdot\cdot\cdot,B_{X_k},
\cdot\cdot\cdot,B_{X_l}) + \Sigma_{E_j \text{not\ appearing\ as\ a\ divisor\
on\ any\ of\ }X_k}\epsilon E_j$$
we have
$$K_{W'} + B_{W'} = (p_k \circ \sigma)^*(K_{X_k} + B_{X_k}) + R_{p_k \circ
\sigma}$$
with
$$\text{supp}\ R_{p_k \circ \sigma} = \text{supp}\ E(p_k \circ \sigma).$$
Hence by Lemma 3.6 (i) each $(X_k,B_{X_k})$ is an end result of $K + B$-MMP
over
$X_k$ starting from $(W',B_{W'})$.  This completes the proof.

\vskip.2in

\proclaim{Theorem 4.5} Let $(X,B_X)$ and $(X',B_{X'})$ be log minimal models
with only ${\Bbb Q}$-factorial wklt singularities in dimension 3.  Suppose they
are Sarkisov related.  Then they are connected by a sequence of log flops
$$(X,B_X) \dashrightarrow (X_1,B_{X_1}) \dashrightarrow \dot\cdot\cdot
(X_k,B_{X_k}) \dashrightarrow \cdot\cdot\cdot (X',B_{X'})$$
where all the $(X_k,B_{X_k})$ are Sarkisov related.

\endproclaim

\demo{proof}\enddemo

An easy application of the Negativity Lemma shows (cf.Koll\'ar[15],Lemma 4.3)
that $(X,B_X)$ and $(X',B_{X'})$ are isomorphic in codimension 1, since they
are Sarkisov related.  Let ${\Cal H}_{X'}$ be a very ample divisor on $X'$ and
${\Cal H}_X$ its strict transform on $X$.  We take a log pair $(W,B_W)$
dominating both log minimal models by birational morphisms $p:(W,B_W)
\rightarrow (X,B_X)$ and $q:(W,B_W) \rightarrow (X',B_{X'})$ as in (iii) of
Proposition 3.5, whose existence is guaranteed by Proposition 4.4.  We claim
that $K_X + B_X + \eta{\Cal H}_X$ is wklt for $0 < \eta << 1$, since
$$Bs({\Cal H}_X) \subset p(R_p)$$
where
$$K_W + B_W = p^*(K_X + B_X) + R_p.$$

If $K_X + B_X + \eta{\Cal H}_X$ is nef, then by Log Abundance (cf.KeMaMc[12])
it is semiample.  Since it is the strict transform of an ample divisor $K_{X'}
+ B_{X'} + \eta{\Cal H}_{X'}$ and since $(X,B_X)$ and $(X',B_{X'})$ are
isomorphic in codimension one and ${\Bbb Q}$-factorial, this implies
$(X,B_X)
\overset\sim\to\rightarrow (X',B_{X'})$.

If $K_X + B_X + \eta{\Cal H}_X$ is not nef, then there is an $K_X + B_X +
\eta{\Cal
H}_X$-negative extremal ray, which must be $K_X + B_X$-trivial and of flopping
type (cf.Koll\'ar[15],Lemma 4.4).  We flop this extremal ray to get
another log minimal model $(X,B_X) \dashrightarrow (X_1,B_{X_1})$ with only
${\Bbb Q}$-factorial wklt singularities.  By construction it is easy to see
that $(X,B_X),(X_1,B_{X_1})$ and $(X',B_{X'})$ are all Sarkisov related and
that $K_{X_1} + B_{X_1} + \eta{\Cal H}_{X_1}$ is wklt.  We proceed
inductively and this procedure has to come to an end, since any sequence of
log flops has to terminate (cf.Shokurov[31]Koll\'ar et al[16]).  Thus
we obtain the desired connecting sequence of log flops between $(X,B_X)$ and
$(X',B_{X'})$.

\vskip.1in

We go back to the discussion of the log Sarkisov program with wklt
singularities.

In the following we prove the well-behavior of the Sarkisov relation for log
Mori fiber spaces with wklt singularities in dimension 2.

\proclaim{Lemma 4.6 = Proposition 3.5 for log Mori fiber spaces with wklt
singularities in dimension 2} Let
$$(X_0,B_{X_o}), (X_1,B_{X_1}), \cdot\cdot\cdot,(X_k,B_{X_k}),
\cdot\cdot\cdot,(X_l,B_{X_l})$$
be log Mori fiber spaces with only ${\Bbb Q}$-factorial wklt singularities in
dimension 2.  Then (i) (ii) and (iii) as in Proposition 3.5 (replacing the
assumption of klt singularities with that of wklt singularities and allowing
the
possibility $\epsilon = 1$) are equivalent.
\endproclaim

\demo{Proof}\enddemo

Again we only have to show the implication $(i) \Rightarrow (iii)$.

Take a log pair $(W,B_W)$ as in (i) and let
$$p_k:(W,B_W) \dashrightarrow (X_k,B_{X_k})$$
be a birational map which is a $K + B$-MMP over $Spec k$.  First observe that
in dimension 2 all the $p_k$ are birational morphisms and that
$$\cup_k {p_k^{-1}}_*(B_{X_k}) \subset B_W$$
is an S.N.C. divisor.
We take a blowup $\sigma:W' \rightarrow W$ with centers over
$$\cup_{E_l \text{not\ appearing\ as\ a\ divisor\ on\ any\ of\ }X_k}E_l$$
until
$$\cup_{E_m \text{not\ appearing\ as\ a\ divisor\ on\ any\ of\ }X_k}E_m \cup
\cup_k {p_k^{-1}}_*(B_{X_k})$$
is an S.N.C. divisor.  Then by setting
$$B_{W'} = D_W(B_{X_0}, \cdot\cdot\cdot,B_{X_k}, \cdot\cdot\cdot,B_{X_l}) +
\Sigma_{E_j \text{not\ appearing\ as\ a\ divisor\ on\ any\ of\ }X_k}E_j,$$
we have
$$K_{W'} + B_{W'} = \sigma^*(K_W + B_W) + R_{\sigma}$$
where $R_{\sigma}$ is an effective ramification divisor (whose support may not
coincide with the exceptional locus $E(\sigma)$).

On the other hand, since for each $k$
$$K_W + B_W = {p_k}^*(K_{X_k} + B_{X_k}) + R_{p_k}$$
where the effective ramification divisor $R_{p_k}$ has the support which
coincides with that of the exceptional locus $E(p_k)$, and since the blowup has
all the centers over
$$\cup_{E_l \text{not\ appearing\ as\ a\ divisor\ on\ any\ of\ }X_k}E_l \subset
\cap_i E(p_i),$$
we conclude
$$K_{W'} + B_{W'} = (p_k \circ \sigma)^*(K_{X_k} + B_{X_k}) + R_{p_k \circ
\sigma}$$
with
$$\text{supp}\ R_{p_k \circ \sigma} = \text{supp}\ E(p_k \circ \sigma).$$
Thus by Lemma 3.6 (i) each $(X_k,B_{X_k})$ is an end result of $K + B$-MMP over
$X_k$ starting from $(W',B_{W'})$.  This completes the proof.

\vskip.1in

In dimension 2, any MMP (log or genuine) is a succession of contractions of
divisors without any flip and thus the resulting Mori fiber space is dominated
by the starting variety through a birational morphism.  Also any $K + B$-MMP
over $Spec\ k$ is a process of $K + B$-MMP over any variety $T$ which is
dominated by relevant log pairs.  These easy observations
unique to dimension 2 make the flowchart for the log Sarkisov with wklt
singularities rather straightforward in dimension 2, compared to higher
dimensional case, where the Sarkisov relation seems more subtle.

\vskip.1in

The log Sarkisov degree $(\mu_k,\lambda_{\epsilon k},e_{\epsilon k})$ is
defined
in the following way.

The quasi-effective threshold $\mu_k$ is as before defined to be the positive
rational number s.t.
$$\mu_k(K_{X_k} + B_{X_k}) + {\Cal H}_{X_k} \equiv 0 \text{\ over\ }S_k.$$

$\lambda_{\epsilon k}$ is defined in the exactly same way as in the case with
klt singularities setting $\epsilon = 1$.

We pay extra attention to how we define $e_{\epsilon k}$.  If we try to
define it in the same way as in the case with klt singularities setting
$\epsilon = 1$, then it would not be well defined, since we may have
infinitely many crepant divisors.  This is one of the difficulties one has to
face once we hit the critical value $\epsilon = 1$ creating the neccessity to
deal with the log canonical locus.

In the flowchart below we show all the divisorial blow ups and intermediate
log Mori fiber spaces are dominated by $(W,B_W)$ that we fix from the
beginning as above satisfying the conditions in (iii) in Proposition 3.5 .

We define $e_{\epsilon k}$ of the intermediate log Mori fiber space that
appear in the due course of the flowchart to be the number of $K + B +
\frac{1}{\lambda_{\epsilon k}}{\Cal H}$-crepant divisors $\bold{ON\ W}$.

\vskip.1in

$$\bold{Flowchart\ for\ Log\ Sarkisov\ Program\ in\ dimenion\ 2}$$
$$\bold{with\ WKLT\ Singularities\ and\ its\ Termination}$$

\vskip.1in

$$\boxed{\text{Case}:\lambda_{\epsilon k} \leq \mu_k}$$

The flowchart for the log Sarkisov program with WKLT singularities in this case
goes
parallel to the one for the genuine Sarkisov program.
After untwisting the birational map by a link of type (III) or (IV), the
quasi-effective threshold strictly decreases.

Moreover, since $\phi_{k+1}:(X_{k+1},B_{X_{k+1}}) \rightarrow S_{k+1}$ is
obtained as an end result of $K + B + \frac{1}{\mu_k}{\Cal H}$-MMP (over
$T$), it follows immediately that $(X_{k+1},B_{X_{k+1}})$ is an end result of
$K + B$-MMP starting from $(W,B_W)$ (over $X_{k+1}$) and thus dominated by
$(W,B_W)$ through a birational morphism.

The claim that there is no infinite number of untwisting (successiove or
unsuccessive) by the links unedr the case $\lambda_{\epsilon k} \leq \mu_k$
can be shown similarly, proving the discreteness of the quasi-effective
thresholds noting that a general fiber of $\phi_k$ is ${\Bbb P}^1$.

\vskip.1in

$$\boxed{\text{Case}:\lambda_{\epsilon k} > \mu_k}$$

We take a maximal divisorial blow up constructed starting from $(W,B_W)$,
which dominates $(X_k,B_{X_k})$ by inductive assumption.

Then just as in the genuine Sarkisov program after untwisting the
birational map by a link of type (II) or (I), the quasi-effective threshold
does not increase
$$\mu_{k+1} \leq \mu_k$$
with equality holding only if
$$\align
\text{either\ }\dim S_{k+1} &> \dim S_k\\
\text{or\ }\dim S_{k+1} &= \dim S_k \text{\ and\ }\psi_k \text{\
is\ square}.\\
\endalign$$
Also it follows similarly that
$$\lambda_{\epsilon, k+1} \leq \lambda_{\epsilon k}.$$
Note that $\phi_{k+1}:(X_{k+1},B_{X_{k+1}}) \rightarrow S_{k+1}$
 is an end result of $K + B + \frac{1}{\lambda_{\epsilon k}}{\Cal H}$-MMP
(over $S_k$) starting from the maximal divisorial blowup $(Z_k, B_{Z_k})$,
which in turn is an end result of $K + B + \frac{1}{\lambda_{\epsilon k}}{\Cal
H}$-MMP possibly followed by $K + B$-MMP (over $X_k$) starting from
$(W,B_W)$.  Therefore, it is easy to see that $(X_{k+1},B_{X_{k+1}})$ is an
end result of $K + B$-MMP (over $X_{k+1}$) starting from $(W,B_W)$.  Thus
$e_{\epsilon, k+1}$ is well -defined and in the above inequality
$$\text{if\ }\lambda_{\epsilon, k+1} = \lambda_{\epsilon k} \text{\ then\
}e_{\epsilon, k+1} \leq e_{\epsilon k} - 1 < e_{\epsilon k}.$$

This immediately proves the claim that there is no infinite (successive)
sequence of untwisting by the links under the case $\lambda_{\epsilon k} >
\mu_k$ with stationary $\lambda_{\epsilon k}$.

Now we take a closer look at the proof of the claim that there is no infinite
(successive) sequence of untwisting by the links under the case
$\lambda_{\epsilon k} > \mu_k$ with stationary quasi-effective threshold
$\mu_k$.

The proof goes parallel for Steps 1 and 2 (cf. Claim 2.2) replacing $K$ with $K
+ B$.  Step 3 becomes meaningless and irrelevant in the case with wklt
singularities and we disregard it, i.e., we don't use Step 3 in our argument
below.  Instead we conclude the argument as follows.  First we remark that the
valuations of $k(X)$ corresponding to the unique exceptional divisors $E_k$ of
the maximal divisorial blowups are all distinct.  Moreover, all the $E_k$ are
divisors on $W$.  On the other hand,
$$a(E_k,X_1,B_{X_1} + \alpha {\Cal H}_{X_1}) \leq a(E_k,X_k,B_{X_k} + \alpha
{\Cal H}_{X_k}) < 0,$$
but there are only finitely many divisors on $W$ with negative discrepancies
w.r.t. $B_{X_1} + \alpha {\Cal H}_{X_1}$, a contradiction!

\vskip.1in

We finish the proof of termination by showing the claim that there is no
infinite (successive) sequence of untwisting by the links under the case
$\lambda_{\epsilon k} > \mu_k$ with nonstationary quasi-effective threshold.

For the case where $\dim S_{k_0} \geq 1$ for some $k_0$ (and thus for
$\forall k \geq k_0$), we show the discreteness of the quasi-effective
thresholds again noting that a general fiber of $\phi_k (k \geq k_0)$ is
${\Bbb P}^1$.

For the case $\forall k, \dim S_k = 0$, we note that the $X_k$ are all
log Del Pezzo surfaces (normal projective surfaces with only quotient
singularities having ample anti-canonical divisors) which are dominated by one
fixed nonsingular projective surface $W$.  Therefore, it is easy to see that
the $X_k$ belong to a bounded family, from which fact the discreteness of the
quasi-effective thresholds follows just as before.

This completes the discussion of the flowchart and its termination for the log
Sarkisov program with wklt singularities in dimension 2.

\vskip.2in

Fianlly we discuss briefly the problems we face when we try to establish the
log Sarkisov program with wklt singularities in higher dimension.

\proclaim{Problem 1}\endproclaim

Does the Sarkisov relation behave well with wklt singularities, i.e., do we
have
the equivalence of (i) (ii) and (iii) in proposition 3.5 replacing klt
singularities with wklt singularities and allowing $\epsilon = 1$?

\vskip.1in

\proclaim{Problem 2}\endproclaim

If the answer to Problem 1 is affirmative, then we can construct a maximal
divisorial blowup of $(X_k,B_{X_k})$ with respect to ${\Cal H}_{X_k}$ in the
case $\lambda_{\epsilon k} \leq \mu_k$ from a good log pair $(W,B_W)$ as in
(iii) of Proposition 3.5.  After $K + B + \frac{1}{\lambda_{\epsilon k}}{\Cal
H}$-MMP we reach $(X_{k+1},B_{X_{k+1}})$.  In the case $\lambda_{\epsilon k}
\leq \mu_k$, after $K + B + \frac{1}{\mu_k}{\Cal H}$-MMP we reach
$(X_{k+1},B_{X_{k+1}})$.

Show in both cases that the $(X_m,B_{X_m})\ m = 0,
1, \cdot\cdot\cdot, k, k+1$ are all Sarkisov related establishing the
inductive procedure.

We note that the N\"other-Fano criterion remains valid as long as we know
$(X,B_X)$, $(X',B_{X'})$ and $(X_k,B_{X_k})$ are Sarkisov related.  (This is
not a trivial remark as at one point of the proof (See Corti[4],Theorem 4.2.)
the positivity of
some coefficient in the ramification divisor does not follow without the
assumption of being
Sarkisov related in the case of wklt singularities.)

\vskip.1in

\proclaim{Problem 3}\endproclaim

Show the discreteness of the quasi-effective thresholds under the case
$\lambda_{\epsilon k} \leq \mu_k$, which follows from the
(conjectural) boundedness of the fibers of the $\phi_k$.

We remark that this is not a straight
consequence of $S_d(Global)$ ($d \leq n - 1$ in dimension $n$), since we do not
know that the number of components in the boundary or the coefficient
$\frac{1}{\mu_k}$ to be bounded.

\vskip.1in

\proclaim{Problem 4}\endproclaim

Show that there is no infinite (successive) sequence of untwisting by the
links under the case $\lambda_{\epsilon k} > \mu_k$ in the following manner:

\ \ i) Show that $\lambda_{\epsilon k}$ cannot be stationary, by adopting an
appropriate definition of $e_{\epsilon k}$ as demonstrated in the case
of dimension 2.  This should be relatively easy.

\ \ ii) Show that $\mu_k$ cannot be stationary.  Steps 1 and 2 of Calim 2.2
go without change, while Step 3 is irrelevant.  We should conclude the
argument by looking at $S_d(Global)\ d \leq n - 1$ on the exceptional divisors
of
the maximal divisorial blowups that appear with coefficient 1.  (This line of
argument was suggested to us by A. Corti.  In fact we could argue this way in
the proof of termination in dimension $n = 2$, though it becomes substantially
lengthier than the one we give.)

\proclaim{Problem 5}\endproclaim

Finally show the discreteness of the quasi-effective thresholds under the case
$\lambda_{\epsilon k} > \mu_k$, which follows again from the boundedness of the
fibers of the
$\phi_k$.  This seems to be the most difficult part.

\newpage

\topmatter
\title References\endtitle
\address Department of Mathematics, Brandeis University, Boston,
Massachusetts 02254-9110\endaddress  \endtopmatter

\widestnumber\no{10}

\ref \no 1 \by V. Alexeev
\paper Two two-dimensional terminations \yr 1993  \vol 69 \jour Duke. Math.
Journal \endref

\ref \no 2 \bysame
\paper Boundedness of $K^2$ for log surfaces \yr 1994   \jour (preprint)
\endref

\ref \no 3 \by A. Borisov
\paper Boundedness theorem for Fano log-threefolds \yr 1994 \jour (preprint)
\endref

\ref \no 4 \by A. Corti
\paper Factorizing birational maps of threefolds after Sarkisov \yr 1993 \jour
(preprint) \endref

\ref \no 5 \by  Y. Kawamata \pages 93--163
\paper Crepant blowing-ups of three dimensional canonical
singularities and its application to degenerations of
surfaces  \yr 1988 \vol 127 \jour Ann. of Math. \endref

\ref \no 6 \bysame \pages 439-445
\paper Boundedness of ${\Bbb Q}$-Fano threefolds \yr 1989 \vol 131
\jour Proc. Int. Conf. Algebra, Contemp. Math. \endref

\ref \no 7 \bysame \pages 229-246
\paper Abundance theorem for minimal threefolds \yr 1991 \vol 108
 \jour Inv. Math. \endref

\ref \no 8 \bysame \pages 609-611
\paper On the length of an extremal rational curve \yr 1991 \vol 105
\jour Invent. math. \endref

\ref \no 9 \bysame \pages 653--660
\paper Termination of log flips for algebraic 3-folds \yr 1992 \vol 3
\jour Internat. J. Math. \endref

\ref \no 10 \by Y.Kawamata-K.Matsuki \pages 595--598 \paper
The number of minimal models for 3-folds of general type is finite. \yr 1987
\vol
276 \jour Math. Ann. \endref

\ref \no 11 \by Y.Kawamata-K.Matsuda-K.Matsuki \pages 283--360 \paper
Introduction to minimal model problem \yr 1985 \vol 10 \jour Adv. Stud. Pure
Math., Alg. Geom., Sendai, T. Oda ed. \endref

\ref \no 12 \by S. Keel-K. Matsuki-J. McKernan \pages 99-119
\paper Log abundance theorem for threefolds \yr 1994 \vol 75
\jour Duke Math. Journal \endref

\ref \no 13 \by S. Keel-J. McKernan
\paper Rational curves on log Del Pezzo surfaces \yr 1994
\jour (preprint) \endref

\ref \no 14 \by K.Kodaira \pages 563--626 \paper On
compact analytic surfaces; II \yr 1963 \vol 77 \jour Ann.
of Math. \endref

\ref \no 15 \by J. Koll\'ar \pages 15--36 \paper
Flops \yr 1989 \vol 113 \jour Nagoya Math. J. \endref

\ref \no 16 \by Koll\'ar et al
\paper Flips and abundance for algebraic threefolds, Summer Seminar Note at the
University of Utah
\yr 1992
\vol 211 \jour Ast\'eridque, J. Koll\'ar edit. \endref

\ref \no 17 \by K. Matsuki
\paper A note on the Sarkisov's program \yr 1992
\jour preprint (unpublished) \endref

\ref \no 18 \bysame
\paper Weyl groups and birational transformations among minimal models \yr 1993
\jour To appear in AMS Memoirs \endref

\ref \no 19 \by Y. Miyaoka \pages 449--476 \paper The chern
classes and Kodaira dimension of a minimal variety \vol 10 \jour Adv.
Stud. Pure Math.\endref

\ref \no 20 \bysame \pages 325--332 \paper On the
Kodaira dimension of minimal threefolds \vol 281 \jour Math. Ann.\endref

\ref \no 21 \bysame \pages 203--220 \paper Abundance Conjecture
for 3-folds: case $\nu = 1$ \vol 68 \jour Compositio Math.\endref

\ref \no 22 \by Y. Miyaoka-S. Mori \pages 65-69 \vol 124
\paper A numerical criterion for uniruledness \jour Ann. of Math. \yr 1986
\endref

\ref \no 23 \by S.Mori \pages 133--176 \paper Threefolds whose canonical
bundles are not numerically effective \yr 1982 \vol 116 \jour Ann. of Math.
\endref

\ref \no 24 \bysame \pages 269--331
\paper Classification of Higher-Dimensional Varieties \yr 1985  \vol 46-Part1
\jour Proc. Sympos. Pure Math. (Summer Research Inst., Bowdoin 1985) \endref

\ref \no 25 \bysame \pages 117--253 \paper Flip theorem and
the existence of minimal models for 3-folds \yr 1988 \vol 1 \jour J. Amer.
Math. Soc. \endref

\ref \no 26 \by V. V. Nikulin \pages 657-675
\paper Del pezzo surfaces with log terminal singularities III \yr 1990 \vol 35
\jour Math. USSR Izv. \endref

\ref \no 27 \by M. Reid \pages 273--310 \paper Canonical threefolds \yr
1980  \jour in Ge\'om\'etrie Alg\'ebrique Angers 1979, A. Beauville ed.,
Sijthoff and Noordhoff \endref

\ref \no 28 \bysame \pages 131--180 \paper Minimal models of
canonical 3-folds \yr 1983 \vol 1 \jour Adv. Stud. in Pure Math. \endref

\ref \no 29 \bysame
\paper Birational geometry of 3-folds according to Sarkisov \yr 1991
\jour (preprint) \endref

\ref \no 30 \by V. G. Sarkisov
\paper Birational maps of standard ${\Bbb Q}$-Fano fiberings \yr 1989
\jour I. V. Kurchatov Institute Atomic Energy preprint \endref

\ref \no 31 \by V. V. Shokurov \pages 105-203
\paper 3-fold log flips \yr 1992 \vol 56
\jour Math. USSR Izv. \endref

\ref \no 32 \by E. Szab\'o \yr 1993
\jour Ph. D. Thesis (University of Utah) \endref

\enddocument